\documentclass[times,doublespace]{article}

\newif\ifInputs

\Inputstrue 

\usepackage{moreverb}

\usepackage{tikz,pgfplots,filecontents}
\usepackage{tkz-fct}
\usepackage{wrapfig}
\usepackage{adjustbox}
\usepackage{amsfonts,amssymb,amsbsy,amsmath,amsthm}
\usepackage{latexsym,graphicx,color}

\newcommand{\half}[1]{\frac{#1}{2}}
\newcommand{\bigO}{\mathcal{O}}

\newcommand{\xx}{{\mathbf{x}}}
\newcommand{\yy}{{\mathbf{y}}}

\newcommand{\ff}{{\mathbf{f}}}
\newcommand{\CC}{{\mathbb{C}}}
\newcommand{\RR}{{\mathbb{R}}}

\newcommand{\DD}{{\mathcal{D}}}

\newcommand{\pderiv}[2]{\frac{\partial #1}{\partial #2}}

\renewcommand{\ss}{{\mathbf{s}}}
\newcommand{\nn}{{\mathbf{n}}}
\newcommand{\eeta}{{\boldsymbol\eta}}

\newcommand{\nG}{{m}}

\begin{document}

\title{On preconditioners for the Laplace double-layer in 2D}

\author{Bryan Quaife and George Biros}

%

\maketitle

\begin{abstract}
 The discretization of the double-layer potential integral equation
 for the interior Dirichlet Laplace problem in a domain with smooth
 boundary results in a linear system that has a bounded condition
 number. Thus, the number of iterations required for the convergence
 of a Krylov method is, asymptotically, independent of the
 discretization size $N$. Using the Fast Multipole Method (FMM) to
 accelerate the matrix-vector products, we obtain an optimal
 $\bigO(N)$ solver. In practice, however, when the geometry is
 complicated, the number of Krylov iterations can be quite large---to
 the extend that necessitates the use of preconditioning.

  We summarize the different methodologies that have appeared in the
  literature (single-grid, multigrid, approximate sparse inverses) and
  we propose a new class of preconditioners based on an FMM-based
  spatial decomposition of the double-layer operator.  We present an
  experimental study in which we compare the different approaches and
  we discuss the merits and shortcomings of our approach.  Our method
  can be easily extended to other second-kind integral equations with
  non-oscillatory kernels in two and three dimensions.
\end{abstract}



\section{Introduction}\label{s:intro}
\begin{figure}[htps]
\centering
\begin{minipage}{0.45\textwidth}
\centering
\ifInputs
\begin{tikzpicture}[scale = 0.8]

\begin{axis}[
  width=2in, height=2in,
  scale only axis,
  xmin=-2, xmax=2,
  ymin=-2, ymax=2,
  hide axis
  ]
\addplot[color=red,line width = 1.0pt,solid,domain=0:360,samples=360*2]({(1+0.98*cos(24*\x))*cos(\x)},{(1+0.98*cos(24*\x))*sin(\x)});
\end{axis}

\end{tikzpicture}
\fi
\end{minipage}
\begin{minipage}{0.45\textwidth}
\centering
\begin{tabular}{rcc}
  $N$ & GMRES & \emph{Preco} GMRES \\
  \hline
  $512  $  & 240  & 9  \\
  $1,024$  & 228  & 11  \\
  $2,048$  & 301  & 16  \\
  $4,096$  & 405  & 20  \\
  $8,192$  & 626  & 24  \\
  $16,382$ & 718  & 29 
\end{tabular}
\end{minipage}
\caption{\label{f:introExample} To illustrate the difficulties of
  ``inverting'' the double-layer potential for a complex geometry
  consider the boundary depicted on the left.  The boundary is
  $C^\infty$, but has high derivatives. Resolving this boundary
  requires a large number of points even when using a spectrally
  accurate discretization. (Note that the boundary does not intersect
  itself in the middle as we can see in the right-most plot of
  Figure~\ref{f:geom}.)  On the right, we report the number of
  unpreconditioned and preconditioned GMRES iterations (\emph{Preco}),
  for tolerance equal to $10^{-12}$, as a function of $N$, the number
  of points used to discretize the boundary.  We have developed a
  novel, FMM-based, preconditioner for integral equations (see
  section~\ref{s:schur-low-rank}), which we used to compute the
  results for the preconditioned GMRES. First, notice that despite the
  fact that we solve a second-kind integral equation, the number of
  GMRES iterations increases. This is because we have not resolved
  sufficiently the geometry.  (Eventually the number of iterations
  does become mesh-independent.) Second, by using our FMM-based
  single-grid preconditioner (with 50 points per leaf) we can
  significantly decrease the number of GMRES steps and the overall
  cost of the calculation.  In section~\ref{s:results}, we will see
  that using multigrid on this problem does not reduce the overall
  computational cost, again because the geometry is barely resolved
  even for $N=16,382$.}
\end{figure}

We consider the solution of the interior Laplace equation with
Dirichlet boundary conditions in a complex, simply-connected domain
$\Omega$ with smooth boundary $\Gamma$, and boundary  data $f$:
\begin{equation*} \Delta u=0, \ \mbox{in}\ \Omega, \quad u=f, \
\mbox{on}\  \Gamma.  \end{equation*} Using a double-layer assumption, we
formulate the problem as a boundary integral equation (BIE)
\begin{equation} \label{e:dlp} \eta(\xx) + \int_\Gamma
K(\xx,\yy)\eta(\yy) = 2f(\xx), \quad \xx \in \Gamma \quad
\mbox{or}\quad (I+\DD)\eta = f, \quad \mbox{or~simply}  \quad A\eta=f.
\end{equation} for the unknown double-layer density $\eta(\yy)$. Here
$K(\xx,\yy)$ is the double-layer kernel of the Laplacian, $I$ is the
identity operator, $\DD$ is the double-layer potential operator, and
the constant $2$ has been absorbed in $f$.

The discretization of~\eqref{e:dlp} results in a dense unsymmetric system of
equations. Gaussian elimination requires $\bigO(N^3)$ work where $N$ is the
number of unknowns, which is prohibitively expensive.  Since~\eqref{e:dlp} has
a bounded condition number~\cite{kress-99}, a Krylov method like the
Generalized Minimum Residual Method (GMRES)~\cite{saa:sch} coupled with an FMM
acceleration~\cite{gre:rok} for the matrix-vector multiplication (henceforth,
``{\em matvec}'') results in $\bigO(N)$ work, which is algorithmically optimal.
The constant in this complexity has two factors, one factor captures the FMM
matvec costs, and the second one captures the number of GMRES iterations. Here
we focus on reducing the latter.

Although the number of GMRES iterations eventually becomes independent
of $N$, this number can be quite large depending on the complexity of
the geometry. (This also implies that in the pre-asymptotic regime the
number of iterations can increase quite significantly.) One way to
quantitatively characterize the complexity of the geometry is by the
number of points $N$ required to approximate $D \eta$ accurately
enough for $\eta(\xx) = 1(\xx)$, the constant density function that is
equal to one for all $\xx$.  If $N$ is large, the solution
of~\eqref{e:dlp} will require hundreds or even thousands of GMRES
iterations. Such cost makes the solution of~\eqref{e:dlp}
prohibitively expensive.  The number of iterations becomes even more
pronounced for boundaries with corners.

\subsection{Summary of methodology and contributions}

Our goal is to design preconditioners that can be applied to a variety
of integral equations with non-oscillatory kernels in two- and
three-dimensions, and can built on top of standard FMM codes.  Based on
the FMM hierarchical decomposition, we propose a new preconditioner,
which we term {\em ``FMMSCHUR''}.  To motivate and evaluate the
preconditioner, we review existing methods and we provide preliminary
comparative results. FMMSCHUR costs $\mathcal{O}(N\log N)$ work and its
construction can be done based on the FMM hierarchy independent of the
dimension of the problem. FMMSCHUR is quite effective in reducing the
number of GMRES iterations when complex geometries are involved.

FMMSCHUR's basic idea is to combine a preconditioner for the near-field
interactions in the FMM with a two-level low rank approximation of the
far-field interactions. In particular it uses the following concepts:
Near-far field decomposition, a Schur approximate factorization based
on this factorization, and algorithms for low-rank approximation of a
matrix. 
\begin{itemize}
\item The {\bf Near-far field decomposition} of $D$  is defined
  as
\begin{align*}
  D=D_{0} + \sum_{\ell} D_{\ell}.
\end{align*}
$D_{0}$ represents the near-field interactions. These also known as
    {\em direct} interactions since they are {\em not} approximated in
    the FMM.  $D_{\ell}$ represent different parts of the far field
    (approximated in the FMM). For example, $D_1$ corresponds to the
    first level of the so called $V$-list interactions.  We make this
    decomposition precise in section~\ref{s:fmm}. Most production
    integral equation codes use single-grid preconditioning techniques
    based on approximate factorizations of $A_0=I+D_{0}$.

\item {\bf Schur low-rank}. Given $A_{1}=A_{0}+D_{1}$ and using the
  FMM factorization of $D_{1}$ which can be written as
  $D_{1}=L_{1}M_{1}^T$, we can use the Sherman-Morrison-Woodbury (SMW)
  formula to approximate the inverse of $A_{1}$.  The operators
  $L_{1}$ and $M^T_{1}$ correspond to the downward pass (translations
  and local-to-local) and upward pass (multipole-to-multipole)
  operators in the FMM. Although formally $D_{1}$ is not low rank, it
  is a compact operator can be approximated using a low-rank
  decomposition (classical FMM is based on this low-rank
  property). Based on this observation, we can approximate the inverse
  of $A_{1}$ as follows.  First we approximate $A_{0}^{-1}$. We use an
  incomplete LU decomposition with a drop tolerance of
  $10^{-3}$. Then, to invert $A_{1}$ using the SMW formula we need to
  approximate the inverse of $I + M_{1}^{T}A_{0}^{-1}L_{1}$. For this,
  we will use SMW one more time and we use a low-rank approximation of
  $M_{1}^{T}A_{0}^{-1}L_{1}$. This approximate inverse significantly
  reduces the number of GMRES iterations required to invert $A_{1}$.
  Finally, once this preconditioner of $A_{1}$ is formed, we use a
  rank $\mathcal{O}(N \log N)$ approximation of $(A - A_{1})$ to
  complete the construction of FMMSCHUR.  This preconditioner can be
  built on top of standard FMM codes, is not restricted to
  two-dimensions, and the overall complexity of applying this
  preconditioner is $\mathcal{O}(N\log N)$.  These ideas are explained
  in detail in section~\ref{s:schur-low-rank}.
\end{itemize}

Our scheme can be further improved as follows. Classical multigrid for
second kind integral equations can be quite effective. In fact for
certain problems it is much more efficient than using a more involved
scheme like FMMSCHUR. In some other problems it pays off to combine
FMMSCHUR with a two-level multigrid.  For a given target accuracy in
the solution, the discretization size $N$ is dictated by $N_f$
required to resolve $f$ and $N_\Gamma$ required to resolve
$\Gamma$. If $N_f\gg N_\Gamma$ multigrid can be used with the coarse
level being the level at which $N\approx N_\Gamma$. If $N_\Gamma \geq
N_f$, multigrid will not help. 

The smoothing scheme in multigrid can
be either the Picard step, namely $\eta =f- D\eta$, or a near-far
field splitting iterative method (see section~\ref{s:preco}). We use a
two-grid scheme combined with a spectral projection-based
discretization method. The coarse grid solve is done with GMRES
preconditioned by FMMSCHUR.

In summary, our method is based on the following formulation and
discretization choices.
\begin{itemize}
\item We use a spectral discrete Galerkin formulation in which the
unknown density $\eta$ is represented by its Fourier coefficients.

\item The boundary integral operator is applied in real space using
Nystr\"{o}m integration with the trapezoid rule.

\item In the multigrid versions of our preconditioners, the restriction
and prolongation operators are spectral projection performed using the
FFT (section~\ref{s:projection}).  We experimented with other
prolongation and restriction operators, but we elected to use spectral
methods so that we can best analyze the quality of the smoothers.

\item Coarsening of $D$ in our multigrid versions is done using two
  alternative methods: projection-based and geometry-based. The former
  is more effective in reducing the number of GMRES iterations but it
  does not reduce the overall cost of the calculation. The latter is
  more efficient computationally as long as the coarse grid operator
  sufficiently resolves the geometry. (The intergrid transfer
  operators are the only component that depends on the particulars of
  the double-layer potential in 2D.)

\item In some of our preconditioning schemes we use low-rank
  approximations of operators. A rank $p$ approximation is computed by
  forming the full SVD $U\Sigma V^{*}$ and then constructing
  $U_{p}\Sigma_{t}V^{*}_{p}$, where $U_{p}$ contains only the $p$
  columns of $U$ and $V_{p}$ the $p$ rows of $V$ that correspond to
  the $p$ largest singular values $\Sigma_{p}$.  This of course is not
  a scalable approach.  Approximate techniques for such decompositions
  can be found~\cite{halko-martinsson-tropp11}, require only
  matrix-vector multiplications and are necessary for obtaining the
  $\bigO(N\log N)$ convergence. For simplicity, we have not
  incorporated such techniques in this initial implementation of
  FMMSCHUR.

\item Our implementation is prototyped in MATLAB so actual run times
are uninformative. To assess the performance of the different
preconditioners we measure their cost in terms of matvecs with $D$.

\end{itemize}

We consider the following preconditioners: multigrid with Picard
smoother, different sparse approximations of $D$, and the FMM-based
decomposition which involves different levels of the FMM operator.  We
define a measure of the complexity of the boundary: the number of
discretization points required to compute $\|D_N 1(\xx)+
\half{1(\xx)}\|/\|1(\xx)\|$ accurately (using either the $L^2$- or the
$L^\infty$-norm).  As mentioned, the effectiveness of the multigrid
depends on how accurately $D_{N}$ approximates $D$ on the coarsest
grid.

{\bf Contributions.} Our first contribution is a numerical study of the
main preconditioning techniques that have appeared in the literature:
(i) We test the classical Picard-multigrid method.  We examine the
effect of approximating the geometry in coarser grids and discuss the
effects of using a projection-based versus a geometry-based coarse grid
operator. {(ii)} We consider sparse approximate preconditioners for $A$
which are based on factorizations of $A_0=I+D_0$. In sparse approximate
preconditioners, one typically introduces two approximations. The first
approximation is that the far field interactions are dropped and only
the near field interactions are used to construct the preconditioner.
The second approximation is that the near field interactions are
factorized using either incomplete factorizations or some variant of a
sparse approximate inverse (SPAI) method. The main trade-offs in
constructing such preconditioners are the computational cost of
constructing and applying the preconditioner, its effectiveness, and
its parallelizability. Here we only consider the issue of
effectiveness.  To assess it, we consider block diagonal, banded,
incomplete, and exact factorizations of the near field. We provide
numerical examples that show that such an approximation reduces the
number of iterations, but the overall speed-up may be less than a
factor of two.

Our second and main contribution is to introduce FMMSCHUR a new
preconditioning scheme that includes the far interactions in the
preconditioner. Our scheme is based on the FMM hierarchy and a
sequence of low-rank approximations combined with the SMW identity. It
resembles the structure of the fast direct
solvers~\cite{greengard-gueyffier-martinsson-rokhlin09}. It uses a
different spatial decomposition which is directly related to the FMM.
It does not lead to an exact factorization but it can be easily built
on top of an FMM scheme.  Experimentally we see a tenfold reduction
on the number of GMRES iterations.  We estimate that the cost of the
preconditioner will be equivalent to two matvecs, so the expected
overall speedup will be a factor of five. We provide numerical
evidence that this preconditioner performs well even when we use a
very low-order approximation of the far field in the preconditioner.

\subsection{Previous Work}
Already in the 60's and 70's, there were attempts to design efficient
iterative solvers for ``inverting'' $I+D$. These efforts continue up to
now with research on preconditioners and direct solvers. There have
been many approaches based on two-grid and multigrid solvers, fixed
point iterations, approximate factorizations of the near field
interactions, analytic preconditioners, and low-order-FMM
self-preconditioning approaches.  Our discussion below is by no means
exhaustive but, we hope, summarizes the main methodologies for
preconditioning the double-layer operator.

Multigrid methods~\cite{hemker-schippers81,hackbusch85,sch,bra:lub}
were one of the first attempts at constructing efficient solvers for
$I+D$.  The main concerns were the choice of the smoothers and the
robustness of the method for Nystr\"{o}m discretizations. The
ingredients for constructing a multigrid scheme are a Picard iteration
as a smoother along with restriction and prolongation using either
injection and Nystr\"{o}m interpolation or a standard
$L^2$-projection.  Further research on multigrid methods extended their
applicability to polygonal domains, high-order piecewise polynomial
Galerkin discretizations, and problems in three dimensions
(see~\cite{atk1994} and the references there in).  The geometries
examined were typically quite simple (spheres, ellipsoids, polygonal
domains with a small number of corners).  We will see that such methods
are indeed optimal, assuming that the coarse grid resolves the geometry
and can be computed inexpensively.  But in many practical cases this is
not the case.

Since the coarse grid solve may not be cheap for complex geometries,
there has been a lot of work on preconditioning integral equations
using single-grid methods, mainly approximate factorizations. It was
quickly recognized that the most natural approach is to approximately
factorize $D_0$, the direct-interaction part of FMM that corresponds to
a sparse matrix.  But exactly factorizing $D_0$ is too expensive using
current technologies.  In~\cite{nabors-white-94}, the authors consider
a preconditioner based on near interactions of the FMM that resembles a
block factorization but retains some off-block interactions; they apply
it successfully to first-kind integral equations. The same
preconditioning idea is used in~\cite{grama-kumar-sameh98}, in which
the authors combine it with a self-preconditioning scheme in which
$I+D$ is preconditioned with a low-order
approximation~\cite{barnes-hut-86} of the operator $D$, combined with
an approximate sparse factorization similar to that
of~\cite{nabors-white-94}. The preconditioner somewhat reduces the
number of GMRES iterations but the overall wall-clock time does not
significantly improve over the unpreconditioned case since inverting
the low-order matrix $I+D$ is also quite expensive.

The ideas of approximate sparse factorizations of the near field
interactions are further explored by Chen~\cite{che} in which banded,
block, and sparse approximate inverse (SPAI) factorizations are
discussed. In the latter, the basic idea of constructing a
preconditioner is to assume a sparsity structure for the preconditioner
and solve for each column of the preconditioner by solving several
small least-squares problems.  It was found that the banded
approximation is not efficient unless the band width approaches the
full width of the matrix. For the sparse approximate inverse, one
difficulty is to find a good sparsity pattern. Band patterns are used
for problems in one and two dimensions (an ellipse) and a cylinder in
three dimensions. More complex geometries were not discussed.
In~\cite{carpentieri-e05}, the SPAI ideas are further pursued in trying
to discover a proper sparsity and truncation (for the least-squares
problem) based on the near interactions constructed by the FMM and
combine them to design a parallel preconditioner.  Their preconditioner
is based on Frobenius-norm minimization with a pattern prescribed in
advance. They solve an extremely challenging problems, 3D Helmholtz
problems in complex geometries (aircraft), and combine the SPAI
preconditioner with the inner-outer iteration of Grama et
al.~\cite{grama-kumar-sameh98}, but it is unclear that the use of the
preconditioner results in significant computational savings. 

Other approaches for preconditioners include incomplete LU (ILU)
factorizations~\cite{lee-zhang-cheng03}, and analytic
preconditioners~\cite{vee:gue:zor:bir} in which the solution of a
regular geometry for which the inverse is known analytically is used to
precondition the operator on an arbitrary simply-connected geometry.  

An alternative to iterative methods is to use a direct method. Of
course this is expensive, so approximate direct methods have been
developed for both integral equations and generic sparse
matrices~\cite{greengard-gueyffier-martinsson-rokhlin09,mar:rok,schmitz-ying12,ho-greengard12}.
These methods are similar to nested dissection and domain decomposition
methods and they rely not on sparsity, but on the fact that if the
boundary of the BIE is decomposed into two pieces, each side can be
solved with only a low-rank update coming from the other side.  This is
applied recursively to form a compressed version of the inverse in
$\bigO(N^{3/2})$ time for highly complex geometries.  For simple
geometries like an ellipsoid, the construction cost becomes $O(N\log
N)$. While this is more expensive than GMRES coupled with the FMM, this
compressed inverse can be applied in $\bigO(N \log N)$ time with a
small constant.  Therefore, for problems that involve multiple
right-hand sides for the same matrix, fast direct solvers are quite
effective.

Finally, in the context of multiply-connected geometries, which we do
not consider here, in~\cite{gre:gre:mcf} the authors use a Picard
preconditioner that is modified to help reduce the additional
iterations that GMRES requires to resolve the extra equations that
arise from the multiply-connectedness of the geometry. They do not
considered this preconditioner in a geometric multigrid setting and
while they consider geometries that are up to 200-ply connected, each
component curve is an ellipse with aspect ratio no greater than three.

In summary, the bulk of the literature can be perhaps roughly
categorized in three classes of solvers. One is to accelerate fixed
point schemes using grid hierarchies. The second is to accelerate GMRES
using some type of sparse factorization based on the near field.
Finally, the recent efforts in fast direct solvers are  promising for
the case of multiple right-hand sides but so far seem to have a
superlinear construction phase.

\subsection{Limitations}

Our methodology is by no means complete. It requires the exact or
approximate factorization of the near field $I+D_0$, whose computation
we do not investigate here. Such a factorization could be computed
using recent methods for fast factorization of sparse
matrices~\cite{schmitz-ying12,martinsson13}.  Unlike fast direct
solvers, our method is approximate and we cannot prove algorithmic
costs, since, for example, we cannot bound the number of GMRES
iterations.

We consider a straightforward problem, the simply-connected interior
Laplace problem in 2D with smooth geometries to isolate the geometry
effects.  Extensions of the proposed methodology in three dimensions,
geometries with corners, and more difficult operators like the
Helmholtz problem will require further work. Several aspects of our
scheme require low-rank properties that do not directly extend to
oscillatory Helmholtz problems.

The coarse grid points are chosen by a simple filtering of the fine
grid.  An adaptive discretization using a smaller number of unknowns
could possibly extend multigrid methods to much more complex
geometries.  Such representations will require more generic (than
spectral) intergrid operators. We have not considered such cases here.

Finally, a more thorough analysis of the preconditioners and smoothers
is also required, but this is difficult to do for arbitrary geometries
especially since we are interested in the actual constants in the
complexity estimates.

\subsection{Organization of the paper} 
In section~\ref{s:formulation} we summarize the discretization, the
structure of the FMM, the decomposition into near and far fields, the
intergrid transfer operators, and the coarse grid operators. In
section~\ref{s:preco} we define the different preconditioning
techniques, and in section~\ref{s:results} we discuss the results.

\section{Formulation}\label{s:formulation}

\begin{table}[htps]
  \centering
  \colorbox{gray!20}{%
    \begin{tabular}{|l|l|}
      \hline
      Notation & Descriptions \\
      \hline
      BIE & Boundary integral equation \\
      FMM & Fast multipole method \\
      GMRES & Generalized minimum residual method\\
      SMW & Sherman-Morrison-Woodbury identity\\
      SVD & Singular value decomposition\\
      \hline 
      $\DD$ & Double-layer potential integral operator \\
      $\kappa$ & Curvature \\
      $\nn$ & Unit outward normal \\
      $K(\xx,\yy)$ & Kernel of the double-layer potential \\ 
      $N$ & Number of points on the finest grid \\
      $N_{\min}$ & Number of points on the coarsest grid \\
      $D_{N}$ & $N$-point discretization of $\DD$ \\
      $m$ & Number of GMRES iterations \\
      $m_{coarse}$ & Number of GMRES iterations on the coarsest grid \\
      $s$ & Maximum number of points in each leaf \\
      $U$-list & The set of boxes that neighbor a leaf node \\
      $V$-list & The set of children of the neighbours of a box's \\
      & parent that are not neighbours themselves \\
      \hline 
    \end{tabular}
  }
  \caption{Frequently used acronyms and symbols.}
\end{table}

\begin{figure}[htps]
  \centering
  \ifInputs
  $\vcenter{\hbox{
      \scalebox{.65}{\input{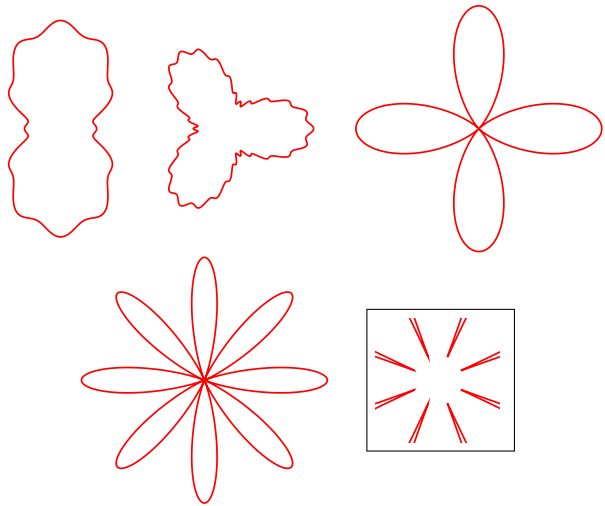}}}}$
  \hspace{2pt}
  $\vcenter{\hbox{
      \scalebox{.65}{\begin{tikzpicture}

\begin{axis}[
  width=1.2in,height=2in,
  scale only axis,
  xmin = -1.0,xmax = 1.6,
  ymin = -1.5,ymax = 1.5,
  hide axis,
  axis equal = true
  ]
\addplot[color=red,line width = 1.0pt,solid,domain=0:360,samples=360*2]({(1+0.5*cos(3*\x)+0.05*cos(30*\x))*cos(\x)},{(1+0.5*cos(3*\x)+0.05*cos(30*\x))*sin(\x)});
\end{axis}

\end{tikzpicture}}}}$
  \hspace{2pt}
  $\vcenter{\hbox{
      \scalebox{.65}{\begin{tikzpicture}

\begin{axis}[
  width=2in, height=2in,
  scale only axis,
  xmin=-2, xmax=2,
  ymin=-2, ymax=2,
  hide axis
  ]
\addplot[color=red,line width = 1.0pt,solid,domain=0:360,samples=360*2]({(1+0.98*cos(4*\x))*cos(\x)},{(1+0.98*cos(4*\x))*sin(\x)});
\end{axis}

\end{tikzpicture}}}}$
  \hspace{2pt}
  $\vcenter{\hbox{
      \scalebox{.65}{\begin{tikzpicture}

\begin{axis}[
  width=2in, height=2in,
  scale only axis,
  xmin=-2, xmax=2,
  ymin=-2, ymax=2,
  hide axis
  ]
\addplot[color=red,line width = 1.0pt,solid,domain=0:360,samples=360*2]({(1+0.98*cos(8*\x))*cos(\x)},{(1+0.98*cos(8*\x))*sin(\x)});
\end{axis}

\end{tikzpicture}}}}$
  \hspace{2pt}
  $\vcenter{\hbox{
      \boxed{\scalebox{.65}{\begin{tikzpicture}

\begin{axis}[
  width=1in,height=1in,
  scale only axis,
  xmin = -0.05,xmax = 0.05,
  ymin = -0.05,ymax = 0.05,
  hide axis
  ]
\addplot[color=red,line width = 1.0pt,solid,domain=0:360,samples=360*2]({(1+0.98*cos(8*\x))*cos(\x)},{(1+0.98*cos(8*\x))*sin(\x)});
\end{axis}

\end{tikzpicture}}}}}$
  \fi
  \caption{\label{f:geom} Four different geometries that we use to
  demonstrate the abilities and limitations of our methods.  From left
  to right, the geometries are as follows.  \emph{Simple:} The
  curvature changes from $-27$ to $17$ at the two rapidly changing
  locations, and with $N=128$ the discretization error is $1.3 \times
  10^{-6}$.  \emph{Moderate:} The curvature ranges from $-188$ to
  $136$, and with $N=256$ the discretization error is $7.2 \times
  10^{-3}.$ \emph{Four-lobed Flower:} The curvature ranges over four
  orders of magnitude and $N=2048$ gives a discretization error of $5.4
  \times 10^{-3}$.  \emph{Eight-lobed Flower:} The curvature ranges
  over six orders of magnitude and $N=4096$ gives a discretization
  error of $7.9 \times 10^{-2}$.  The final plot is a 40 times
  magnification of the center of the eight-lobed flower.}
\end{figure}
Boundary Integral Equations (BIEs) are a reformulation of boundary
value problems that have the advantage of only having to discretize the
boundary of the physical domain.  This reduces the dimension of the
problem by one and allows for complex domains.  Here we formalize the
problem of interest and set notation.  Let $\Omega \subset \RR^{2}$ be
a bounded, simply-connected, domain with smooth boundary $\Gamma$.  The
geometries that we will consider are $\xx(\theta) =
(c r(\theta)\cos(\theta),r(\theta)\sin(\theta))$, $0 \leq \theta < 2\pi,$
where
\begin{align*}
  \emph{Simple}: \quad & r(\theta) = 0.5(\cos^{2}(\theta)+9\sin^{2}(\theta))^{1/2} 
      + 0.07\cos(12\theta),\, c=0.85, \\
  \emph{Moderate}: \quad & r(\theta) = 1 + 0.5\cos(3\theta) + 0.05\cos(30\theta),\, c=1, \\
  \emph{$k$-lobed Flower}: \quad & r(\theta) = 1 + 0.98\cos(k\theta),\, c=1. \\
\end{align*}
Figure~\ref{f:geom} illustrates these geometries that we are considering.
Consider the boundary value problem
\begin{align}
  \Delta u = 0, \ \xx \in \Omega, \quad u = g, \ \xx \in \Gamma.
  \label{e:lap:bvp}
\end{align}
Using the fundamental solution $-\frac{1}{2\pi}\log|\xx|$ of~\eqref{e:lap:bvp},
we seek a solution in the form
\begin{align}
  u(\xx) &= -\frac{1}{2\pi}\int_{\Gamma} \pderiv{}{\nn_{\yy}}\log|\xx - \yy|
  \eta(\yy)ds_{\yy}, \quad \xx \in \Omega,
  \label{e:lap:bie}
\end{align}
where $\nn_{\yy}$ is the unit outward normal to $\Gamma$ at $\yy$.
The ansatz~\eqref{e:lap:bie} is called the double-layer potential and
we write $u(\xx) = \DD[\eta](\xx)$.  The function
\begin{align*}
  K(\xx,\yy) = -\frac{1}{2\pi} \pderiv{}{\nn_{\yy}}\log|\xx - \yy|
\end{align*}
is called the kernel of the integral operator, and its limiting value
is 
\begin{align}
  \lim_{\substack{\yy \in \Gamma \\ \yy \rightarrow \xx_{0}}} 
  K(\xx_{0},\yy) = \frac{1}{4\pi}\kappa(\xx_{0}),
  \label{e:diagonal:term}
\end{align}
where $\kappa(\xx_{0})$ is the curvature at $\xx_{0} \in \Gamma$.  To
satisfy the boundary data of~\eqref{e:lap:bvp}, $\eta$ must satisfy the
Fredholm integral equation of the second kind~\cite{fol,kel,mik}
\begin{align*}
  g(\xx_{0}) = \frac{1}{2}\eta(\xx_{0}) + \DD[\eta](\xx_{0}), \quad \xx_{0} \in \Gamma.
\end{align*}
By defining $f=2g$ and absorbing the $2$ into the double-layer potential, we need to solve 
\begin{align}
  f(\xx_{0}) = \eta(\xx_{0}) + \DD[\eta](\xx_{0}), \quad \xx_{0} \in \Gamma.
  \label{e:fredholm:2kind}
\end{align}
Since the kernel $K$ is bounded, the integral operator $\DD$ is compact
and the Fredholm alternative guarantees the existence and uniqueness of
solutions of~\eqref{e:fredholm:2kind}.

\subsection{Discretization}
To discretize~\eqref{e:fredholm:2kind}, we let $\xx(\theta)$ be a parameterization of $\Gamma$ and discretize $\Gamma$ at $\{\xx(\theta_{j}) \:|\: \theta_{j}=(j-1)2\pi/N,j=1,\ldots,N\}$.  We let $\ff_{j} = \ff(\xx(\theta_{j}))$ and use similar notation for the other variables.  The resulting linear system is
\begin{align}
  \ff_{j} = \eeta_{j} + \sum_{k=1}^{N} K_{jk} \eeta_{k}, \ j=1,\ldots,N, 
  \quad \mbox{or} \quad \ff = \left(I + D_{N}\right)\eeta,
  \label{e:fredholm:discrete}
\end{align}
where $K_{jk} = 2K(\xx_{j},\xx_{k})|\Delta \ss_{k}|$, and $\Delta \ss$
is the Jacobian.  Equation~\eqref{e:diagonal:term} is used for the
diagonal terms $K_{jj}$.  Equation~\eqref{e:fredholm:discrete} can be
interpreted as a Fourier projection method since the $N$-point
trapezoid rule integrates the first $N$ Fourier frequencies exactly.
That is, if we let $F_{N}^{T}$ denote the $N$-point Fourier projection
and let $F_{N}$ be its inverse, then~\eqref{e:fredholm:discrete} is
equivalent to
\begin{align*}
  F_{N}^{T} (I + D_{N}) F_{N} F_{N}^{T}\eeta_{N}= 
      F_{N}^{T} \ff_{N}.
\end{align*}
This motivates using Fourier restriction and prolongation operators which we define in the next section.

Since~\eqref{e:fredholm:2kind} is a second-kind BIE (it would be of the
first kind if the term $\eta(\xx_{0})$ was not present), we are
guaranteed that $\nG$ is bounded independent of $N$.  The number of
required GMRES iterations, $\nG$, depends only on the complexity of
$\Gamma$.  In Table~\ref{t:aspect:ratio}, we consider ellipses with
different aspect ratios, and in Table~\ref{t:number:lobes} we consider
the flower-shaped geometry (see Figure~\ref{f:geom}) with different
numbers of lobes.  We see that the largest curvature has a slight affect
on the number of GMRES iterations but it is the number of high
curvature points that substantially increases $\nG$.  However, given a
fixed geometry, once $N$ is sufficiently large that the geometry is
resolved, mesh independence is observed from this point onwards.

\begin{table}[htps]
  \begin{minipage}{0.4\textwidth}
    \centering
    \begin{tabular}{cc}
      Aspect Ratio & $\nG$ \\
      \hline
      2 & 3 \\
      4 & 4 \\
      8 & 7 \\
      16 & 10 \\
      32 & 13 \\
      64 & 15 \\
      128 & 17 
    \end{tabular}
    \caption{\label{t:aspect:ratio}The number of required GMRES
      iterations to achieve a tolerance of $10^{-12}$ for ellipses of
      varying aspect ratios.  We see that the aspect ratio has a small
      role on the number of GMRES iterations.}
  \end{minipage}
  \hspace*{\fill}
  \begin{minipage}{0.4\textwidth}
    \centering
    \begin{tabular}{cccccccc}
      Number of Lobes & $\nG$ \\
      \hline
      2 & 28 \\
      3 & 44 \\
      4 & 53 \\
      5 & 90 \\
      6 & 91 \\
      7 & 127 \\
      8 & 171
    \end{tabular}
    \caption{\label{t:number:lobes}The number of required GMRES
    iterations to achieve a tolerance of $10^{-12}$ for a flower-shaped
    geometry with a varying number of lobes.  We see that the number of
    points of high curvature has an important role on the number of
    GMRES iterations.}
  \end{minipage}
\end{table}

\subsection{Restriction and Prolongation}\label{s:projection}
We introduce the spectral prolongation $P$ and restriction $R$ operators.
Restriction is done by taking the Fourier transform, truncating the tail of the
spectrum, and then computing the inverse Fourier transform.  Mathematically, if
$\ff \in \CC^{N}$ is represented as
\begin{align*}
  \ff_{j} = \sum_{k=-N/2}^{N/2-1} \hat{f}_{k} e^{ik\theta_{j}}, 
  \quad \theta_{j} = (j-1)\frac{2\pi}{N},\, j=1,\ldots,N,
\end{align*}
then the restriction of $\ff$ is
\begin{align*}
  (R_{N}^{N/2} \ff)_{j} = \sum_{k=-N/4}^{N/4-1} \hat{f}_{k} e^{ik\theta_{j}},
  \quad \theta_{j} = (j-1)\frac{4\pi}{N},\, j=1,\ldots,N/2.
\end{align*}
Similarly, if $\ff \in \CC^{N/2}$ is represented as
\begin{align*}
  \ff_{j} = \sum_{k=-N/4}^{N/4-1} \hat{f}_{k} e^{ik\theta_{j}}, 
  \quad \theta_{j} = (j-1)\frac{4\pi}{N},\, j=1,\ldots,N/2,
\end{align*}
then the prolongation of $\ff$ is
\begin{align*}
  (P_{N/2}^{N} \ff)_{j} = \sum_{k=-N/4}^{N/4-1} \hat{f}_{k} e^{ik\theta_{j}},
  \quad \theta_{j} = (j-1)\frac{2\pi}{N},\, j=1,\ldots,N.
\end{align*}

\subsection{Coarse Grid Operators}\label{s:coarse}
To use a multigrid method, we also need to construct the double-layer potential
$D$ on coarser grids.  We will test two common methods of constructing these
operators. The first method is geometry-based.  It uses $R$ to coarsen
$\Gamma$ and then constructs $D$ using~\eqref{e:fredholm:discrete}.  The second
method is projection-based.  Given the fine grid operator $D_{N} \in \RR^{N
\times N}$, the coarse grid operator for $M < N$ is
\begin{align*}
  D_{M} = R_{N}^{M} D_{N} P_{M}^{N}.
\end{align*}
The projection-based coarse grid operator is more accurate, however, it is not
practical for large problems since the coarse grid operators require
calculations on the finest grid.  Since we are using Fourier based restriction
and prolongation, we always let $N/M$ be a power of $2$.

\subsection{Fast Multipole Method}\label{s:fmm}
\begin{figure}[htps]
  \centering
  \includegraphics[trim = 5.0cm 8.0cm 4.5cm 7.5cm, clip = truee, scale = 0.35]{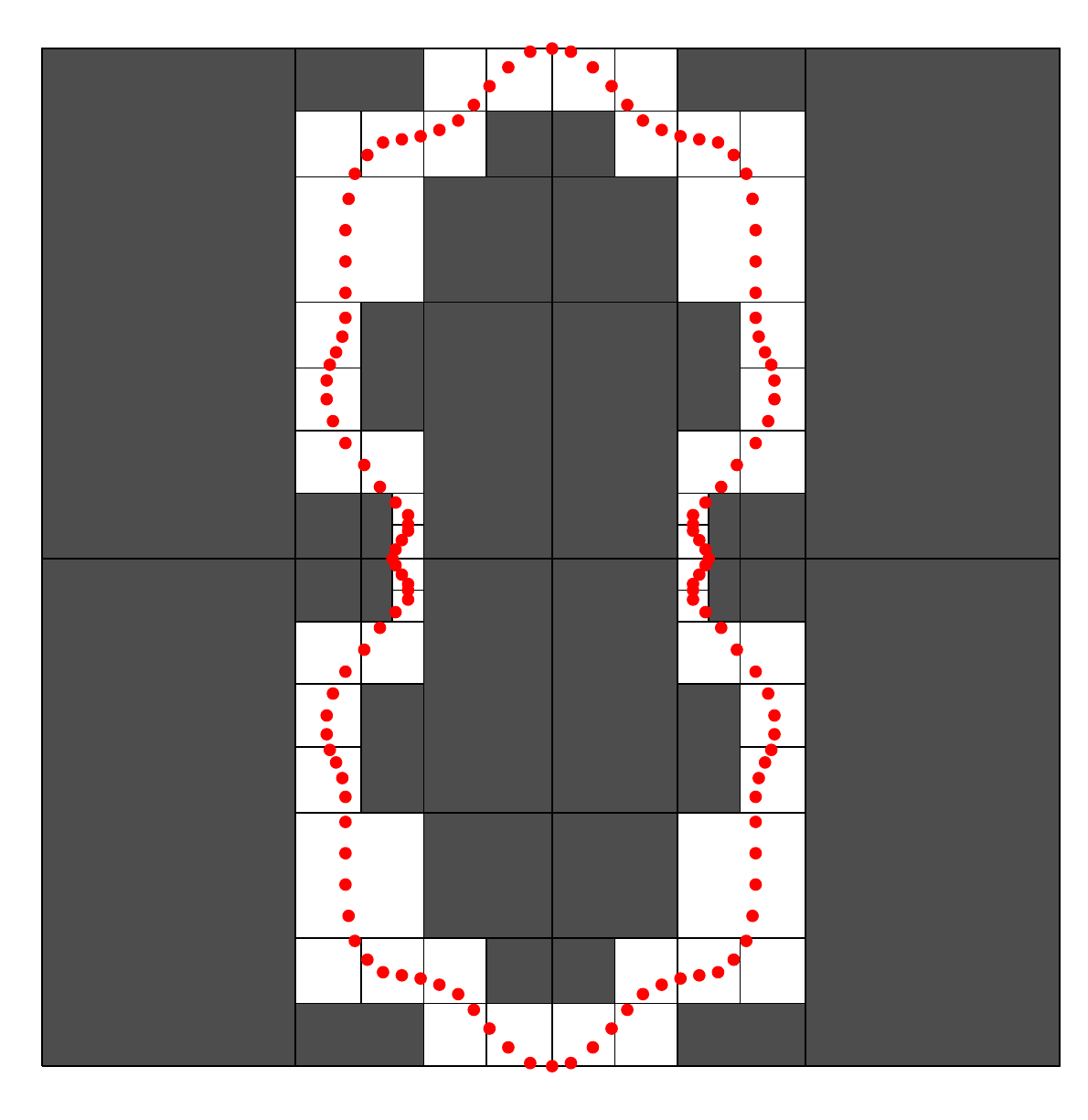}
  \includegraphics[trim = 5.0cm 8.0cm 4.5cm 7.5cm, clip = true, scale = 0.35]{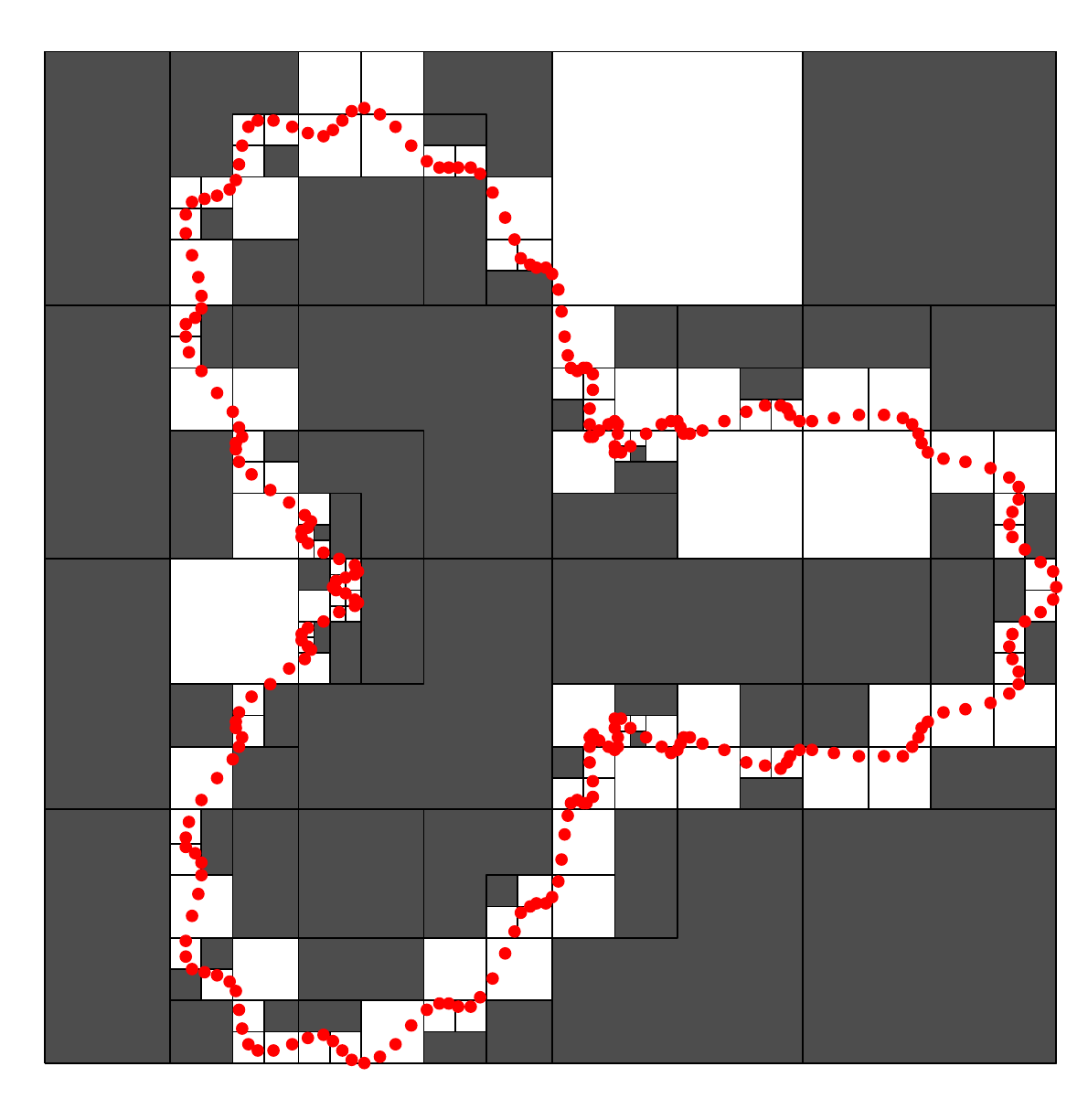}
  \includegraphics[trim = 5.0cm 8.0cm 4.5cm 7.5cm, clip = true, scale = 0.35]{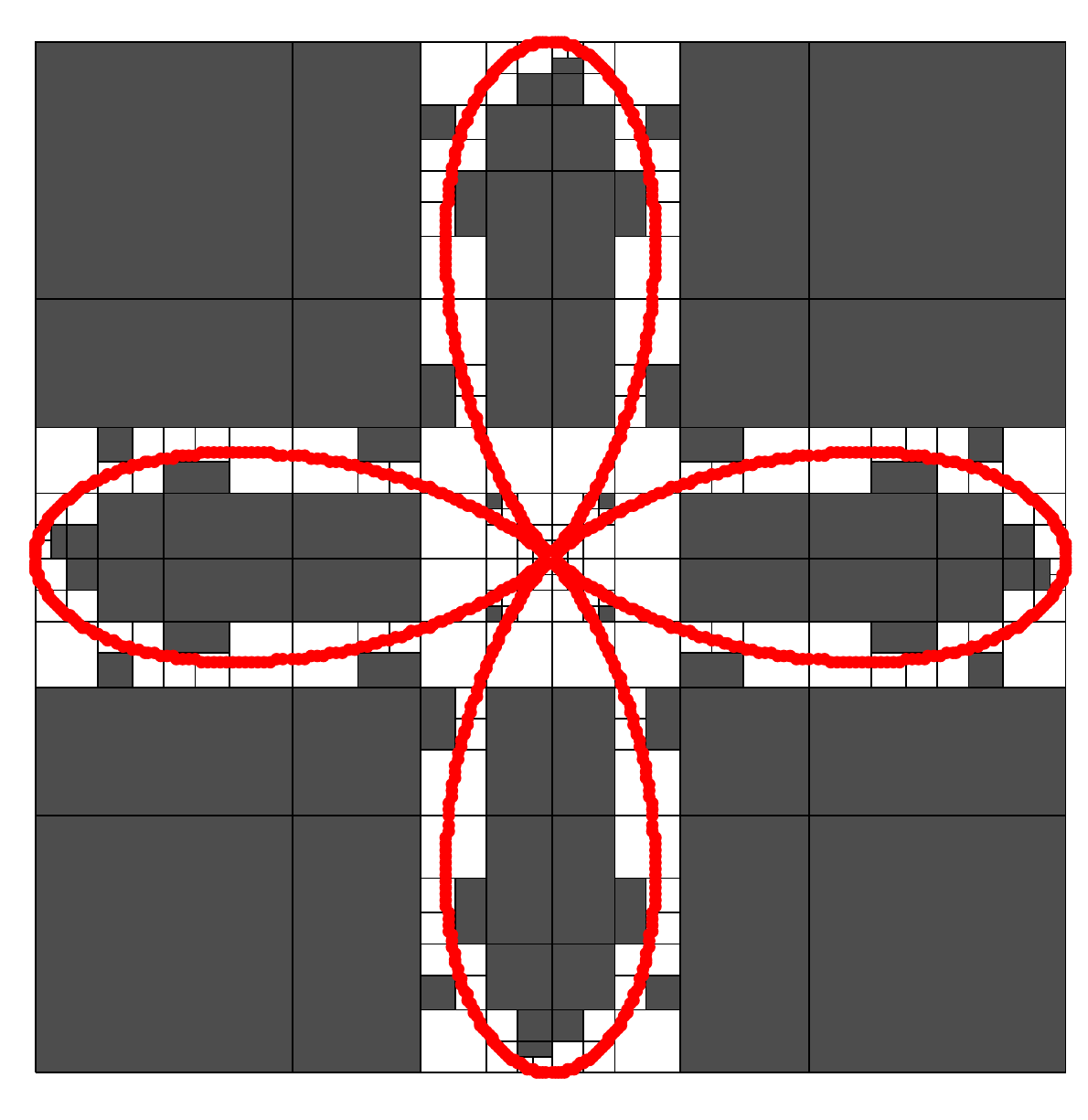}
  \caption{\label{f:quadtrees}The quadtrees for three of the
    geometries we introduced in Figure~\ref{f:geom}.  In the left two
    plots, each leaf contains no more than 4 points.  In the right
    plot, each leaf contains no more than 10 points.  Darkened regions
    contain no points and are removed from the calculation.}
\end{figure}
The Fast Multipole Method (FMM), first introduced by Greengard and
Rokhlin in~\cite{gre:rok}, is a hierarchical method that can apply the
matrix-vector multiplication in~\eqref{e:fredholm:discrete} with
$\bigO(N)$ operations.  We briefly highlight the components of the FMM
that are relevant to our work.  Given a discretization of $\Gamma$, we
form a quadtree so that no leaf contains more than $s$ points.  Our
implementation handles adaptive quadtrees (see
Figure~\ref{f:quadtrees}), but to minimize complex definitions, we only
present the formulation for uniform quadtrees.

\begin{figure}[htps]
  \centering
  \begin{tikzpicture}[xscale=0.4,yscale=0.4]

\draw (0,0)--(0,12);
\draw (0,12)--(12,12);
\draw (12,12)--(12,0);
\draw (0,0)--(12,0);

\draw (4,0)--(4,12);
\draw (8,0)--(8,12);
\draw (0,4)--(12,4);
\draw (0,8)--(12,8);

\draw (2,0)--(2,12);
\draw (6,0)--(6,12);
\draw (10,0)--(10,12);
\draw (0,2)--(12,2);
\draw (0,6)--(12,6);
\draw (0,10)--(12,10);

\draw[ultra thick] (4,4)--(8,4)--(8,8)--(4,8)--cycle;

\node at (5,5) {$b$};
%
\foreach \x in {1,3,5,7,9,11}
{
\node at (\x,1) {$V_{1}$};
\node at (\x,9) {$V_{1}$};
\node at (\x,11) {$V_{1}$};
}

\foreach \y in {3,5,7}
{
\node at (1,\y) {$V_{1}$};
\node at (9,\y) {$V_{1}$};
\node at (11,\y) {$V_{1}$};
}

\node at (3,3) {$U$};
\node at (3,5) {$U$};
\node at (3,7) {$U$};
\node at (5,3) {$U$};
\node at (5,7) {$U$};
\node at (7,3) {$U$};
\node at (7,5) {$U$};
\node at (7,7) {$U$};


\end{tikzpicture}
  \caption{\label{f:uniformQuadtree} The $U$-list and $V_{1}$-list of
  box $b$.  The parent box $P(b)$ is in bold.  The boxes that are not
  illustrated are in $V_{2},V_{3},\ldots$}
\end{figure}
We start by defining two boxes as being neighbours if they share at
least one corner or one edge.  We say that a box is a neighbour of
itself.  Given a leaf $b$ in the quadtree, we let $P(b)$ be its
parent.  Then, we divide the quadtree into multiple regions: the
\emph{$U$-list} and the \emph{$V_{i}$-lists}.  The $U$-list is the set
of neighbors that are adjacent to $b$, and the $V_{1}$-list is the set
of children of $P(b)$'s neighbours that are not neighbors of $b$ (see
Figure~\ref{f:uniformQuadtree}).  We write $U$ for the set of boxes in
the $U$-list, $V_{1}$ for the set of boxes in the $V_{1}$-list, $V_{2}$
for the set of boxes in the $V_{1}$-list of $P(b)$, etc.  Also, we
define $T(B)$ to be the set of points in the set of boxes $B$.

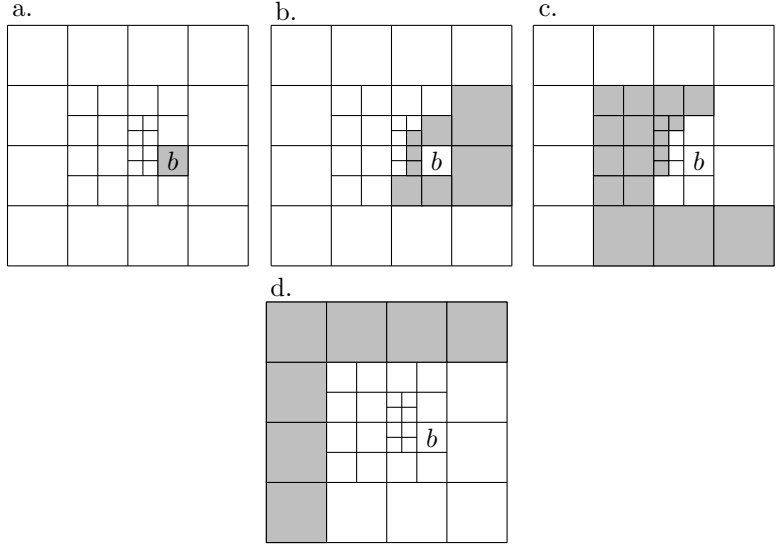
\begin{figure}[htps]
  \centering
  \ifInputs
  \begin{tikzpicture}[xscale=0.2,yscale=0.2]

\draw (0,0)--(0,16);
\draw (0,16)--(16,16);
\draw (16,16)--(16,0);
\draw (0,0)--(16,0);

\draw (4,0)--(4,16);
\draw (8,0)--(8,16);
\draw (12,0)--(12,16);
\draw (0,4)--(16,4);
\draw (0,8)--(16,8);
\draw (0,12)--(16,12);

\draw (6,4)--(6,12);
\draw (10,4)--(10,12);
\draw (4,6)--(12,6);
\draw (4,10)--(12,10);

\draw (9,6)--(9,10);
\draw (8,7)--(10,7);
\draw (8,9)--(10,9);

\draw[fill=gray!50] (10,6)--(12,6)--(12,8)--(10,8)--cycle;

\node at (11,7) {$b$};
\node at (1,17) {a.};


\end{tikzpicture} 
  \begin{tikzpicture}[xscale=0.2,yscale=0.2]

\draw (0,0)--(0,16);
\draw (0,16)--(16,16);
\draw (16,16)--(16,0);
\draw (0,0)--(16,0);

\draw (4,0)--(4,16);
\draw (8,0)--(8,16);
\draw (12,0)--(12,16);
\draw (0,4)--(16,4);
\draw (0,8)--(16,8);
\draw (0,12)--(16,12);

\draw (6,4)--(6,12);
\draw (10,4)--(10,12);
\draw (4,6)--(12,6);
\draw (4,10)--(12,10);

\draw (9,6)--(9,10);
\draw (8,7)--(10,7);
\draw (8,9)--(10,9);

\draw[fill=gray!50] (12,4)--(16,4)--(16,8)--(12,8)--cycle;
\draw[fill=gray!50] (12,8)--(16,8)--(16,12)--(12,12)--cycle;
\draw[fill=gray!50] (10,8)--(12,8)--(12,10)--(10,10)--cycle;
\draw[fill=gray!50] (10,4)--(12,4)--(12,6)--(10,6)--cycle;
\draw[fill=gray!50] (8,4)--(10,4)--(10,6)--(8,6)--cycle;
\draw[fill=gray!50] (9,6)--(10,6)--(10,7)--(9,7)--cycle;
\draw[fill=gray!50] (9,7)--(10,7)--(10,8)--(9,8)--cycle;
\draw[fill=gray!50] (9,8)--(10,8)--(10,9)--(9,9)--cycle;

\node at (11,7) {$b$};
\node at (1,17) {b.};


\end{tikzpicture}
  \begin{tikzpicture}[xscale=0.2,yscale=0.2]

\draw (0,0)--(0,16);
\draw (0,16)--(16,16);
\draw (16,16)--(16,0);
\draw (0,0)--(16,0);

\draw (4,0)--(4,16);
\draw (8,0)--(8,16);
\draw (12,0)--(12,16);
\draw (0,4)--(16,4);
\draw (0,8)--(16,8);
\draw (0,12)--(16,12);

\draw (6,4)--(6,12);
\draw (10,4)--(10,12);
\draw (4,6)--(12,6);
\draw (4,10)--(12,10);

\draw (9,6)--(9,10);
\draw (8,7)--(10,7);
\draw (8,9)--(10,9);

\draw[fill=gray!50] (9,9)--(10,9)--(10,10)--(9,10)--cycle;
\draw[fill=gray!50] (8,6)--(9,6)--(9,7)--(8,7)--cycle;
\draw[fill=gray!50] (8,7)--(9,7)--(9,8)--(8,8)--cycle;
\draw[fill=gray!50] (8,8)--(9,8)--(9,9)--(8,9)--cycle;
\draw[fill=gray!50] (8,9)--(9,9)--(9,10)--(8,10)--cycle;
\draw[fill=gray!50] (8,10)--(10,10)--(10,12)--(8,12)--cycle;
\draw[fill=gray!50] (8,10)--(10,10)--(10,12)--(8,12)--cycle;
\draw[fill=gray!50] (10,10)--(12,10)--(12,12)--(10,12)--cycle;
\draw[fill=gray!50] (4,4)--(6,4)--(6,6)--(4,6)--cycle;
\draw[fill=gray!50] (4,6)--(6,6)--(6,8)--(4,8)--cycle;
\draw[fill=gray!50] (4,8)--(6,8)--(6,10)--(4,10)--cycle;
\draw[fill=gray!50] (4,10)--(6,10)--(6,12)--(4,12)--cycle;
\draw[fill=gray!50] (6,4)--(8,4)--(8,6)--(6,6)--cycle;
\draw[fill=gray!50] (6,6)--(8,6)--(8,8)--(6,8)--cycle;
\draw[fill=gray!50] (6,8)--(8,8)--(8,10)--(6,10)--cycle;
\draw[fill=gray!50] (6,10)--(8,10)--(8,12)--(6,12)--cycle;
\draw[fill=gray!50] (4,0)--(8,0)--(8,4)--(4,4)--cycle;
\draw[fill=gray!50] (8,0)--(12,0)--(12,4)--(8,4)--cycle;
\draw[fill=gray!50] (12,0)--(16,0)--(16,4)--(12,4)--cycle;

\node at (11,7) {$b$};
\node at (1,17) {c.};


\end{tikzpicture}
  \begin{tikzpicture}[xscale=0.2,yscale=0.2]

\draw (0,0)--(0,16);
\draw (0,16)--(16,16);
\draw (16,16)--(16,0);
\draw (0,0)--(16,0);

\draw (4,0)--(4,16);
\draw (8,0)--(8,16);
\draw (12,0)--(12,16);
\draw (0,4)--(16,4);
\draw (0,8)--(16,8);
\draw (0,12)--(16,12);

\draw (6,4)--(6,12);
\draw (10,4)--(10,12);
\draw (4,6)--(12,6);
\draw (4,10)--(12,10);

\draw (9,6)--(9,10);
\draw (8,7)--(10,7);
\draw (8,9)--(10,9);

\draw[fill=gray!50] (0,0)--(4,0)--(4,4)--(0,4)--cycle;
\draw[fill=gray!50] (0,4)--(4,4)--(4,8)--(0,8)--cycle;
\draw[fill=gray!50] (0,8)--(4,8)--(4,12)--(0,12)--cycle;
\draw[fill=gray!50] (0,12)--(4,12)--(4,16)--(0,16)--cycle;
\draw[fill=gray!50] (4,12)--(8,12)--(8,16)--(4,16)--cycle;
\draw[fill=gray!50] (8,12)--(12,12)--(12,16)--(8,16)--cycle;
\draw[fill=gray!50] (12,12)--(16,12)--(16,16)--(12,16)--cycle;

\node at (11,7) {$b$};
\node at (1,17) {d.};


\end{tikzpicture}
  \fi
  \caption{\label{f:mgLayers} For points inside box $b$, the layer
  potential $D$ is approximated with points in a set of neighboring
  boxes.  a. $D_{0}$: Box $b$; b. $D_{0}$: Boxes in $U$; c.  $D_{1}$:
  Boxes in $V_{1}$; d. $D_{2}$: Boxes in $V_{2}$.}
\end{figure}

Given a vector $\eeta$, we decompose the evaluation of $D\eeta$ as
\begin{align*}
  (D\eeta)_{j} = \sum_{k=1}^{N} K_{jk}\eeta_{k} &= \sum_{\xx_{k} \in T(U)} K_{jk} \eeta_{k} 
  &&+ \sum_{\xx_{k} \in T(V_{1})} K_{jk} \eeta_{k} &&+ \sum_{\xx_{k} \in T(V_{2})} K_{jk} \eeta_{k} &&+ \cdots \\
  &= \hspace{22pt}(D_{0}\eeta)_{j} &&+ \hspace{22pt}(D_{1}\eeta)_{j} &&+ \hspace{22pt}(D_{2}\eeta)_{j} &&+ \cdots
\end{align*}
In Figure~\ref{f:mgLayers}, we illustrate the set of points in the
summations for an adaptive tree structure.  The adaptivity changes the
definitions of $V_{1},V_{2},\ldots$ However, to avoid introducing
additional notation, we do not redefine these sets and refer the reader
to Figure~\ref{f:mgLayers} to visualize the different sets.  The
summation due to points in $U$, shown in the left two plots, is denoted
by $D_{0}$.  Here we have made a distinction between $b$ and its other
neighbours.  This will play a role when we introduce our
preconditioners.  The summation due to the points in $V_{1}$, shown in
the third plot, is denoted by $D_{1}$.  The summation due to the points
in $V_{2}$, shown in the final plot, is denoted by $D_{2}$.

The key idea of the FMM is methods for approximating the far field. The
matrices $D_\ell, \ell\geq 1$ can be approximated by a factorization of the
form 
\begin{align*}
  D_\ell =  L_\ell T_\ell M_\ell^T.
\end{align*}
$M_\ell^T$ represents building the multipole moments (or the equivalent
densities for kernel-independent
methods~\cite{yin:bir:zor2004,sun-pitsianis01}), $T_\ell$ represents
translation operators that transform multipole coefficients to
polynomial coefficients, and $L_\ell$ represents the downward traversal
of a tree. For the rest of the paper we absorb $T_\ell$ in $L_\ell$ so
that $D_\ell = L_\ell M_\ell^T$. 

$L_\ell,M_\ell$ can be approximated efficiently in the FMM to result in
an $\bigO(N)$ scheme.  So, $D_\ell$ is an $N\times N$ matrix and can be
approximated very well by $L_\ell M_\ell^T$, where $L_\ell, M_\ell \in
\RR^{N \times m_\ell}$.  For example, if we use $s$ points per box and
we use $p$ multipole moments as we traverse the tree, $m_\ell=p
\frac{N}{s 4^{\ell-1}}$.

\section{Preconditioners}\label{s:preco}
In this section, we introduce several strategies for solving and
preconditioning~\eqref{e:fredholm:discrete}. We focus mainly on
single-grid preconditioners.

\subsection{Picard}
The simplest iterative method for solving~\eqref{e:fredholm:discrete}
is using the Picard iteration
\begin{align}
  \eeta^{\mathrm{new}} = f - D\eeta^{\mathrm{old}}.
  \label{e:picard}
\end{align}
The Picard iteration~\eqref{e:picard} is not convergent since standard
potential theory~\cite{kel} shows that $\DD[1](\xx) = 1$.  Therefore,
as $N \rightarrow \infty$, the spectral radius of $D$ approaches 1
and~\eqref{e:picard} is only stable.  Of course~\eqref{e:picard} is not
being used as an iterative solver.  Instead, it is used as a
smoother~\cite{kress-99,hackbusch85,atk1992,atkinson-97} in a multigrid
scheme.

By coupling this smoother with the previously discussed prolongation,
restriction, and coarse grid operators, a multigrid algorithm is
formed.  Multigrid can then be used as a solver
for~\eqref{e:fredholm:discrete}.  However, it is more efficient to use
it as a preconditioner for GMRES.  Unfortunately, using the Picard
smoother when the coarse grid is not sufficiently resolved results in
a multigrid algorithm that can slow down the convergence of GMRES.

\subsection{Near-Far Field Splitting}\label{s:near}

\begin{figure}
\centering
\fbox{\includegraphics[scale = 0.7,trim = 5.0cm 8.0cm 4.5cm 7.5cm, clip = true, scale = 0.35]{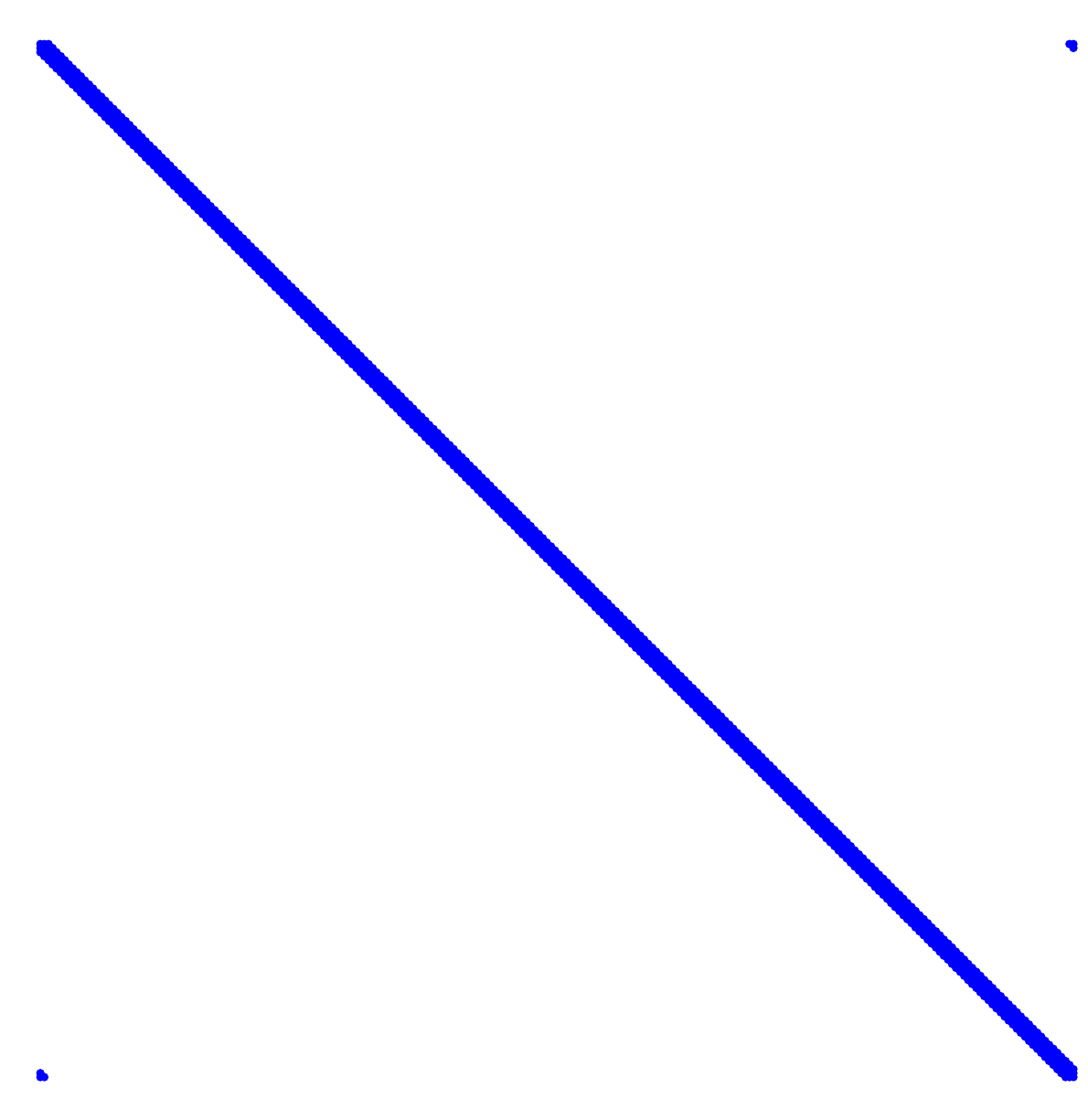}}
\fbox{\includegraphics[scale = 0.7,trim = 5.0cm 8.0cm 4.5cm 7.5cm, clip = true, scale = 0.35]{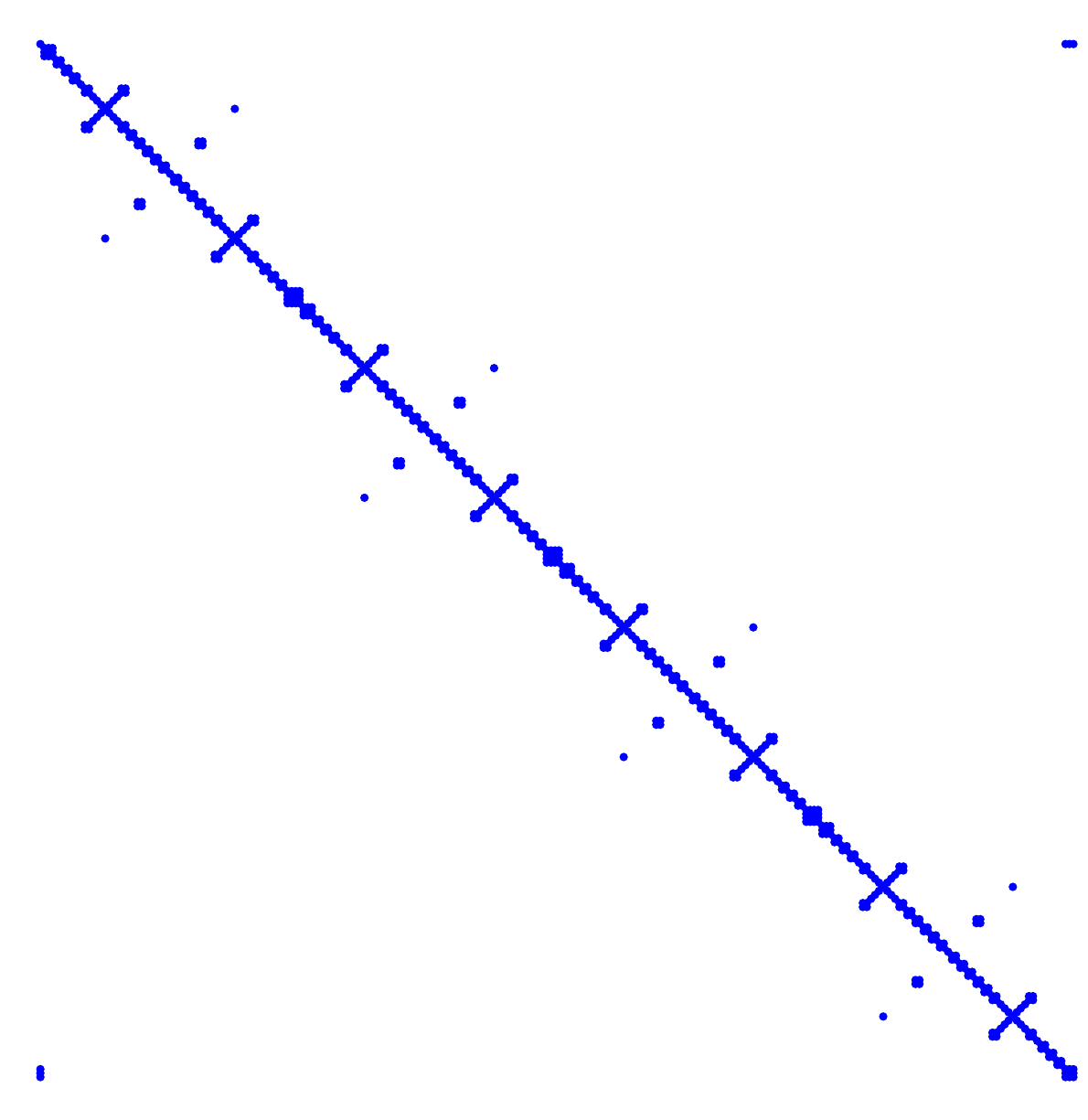}}
\fbox{\includegraphics[scale = 0.7,trim = 5.0cm 8.0cm 4.5cm 7.5cm, clip = true, scale = 0.35]{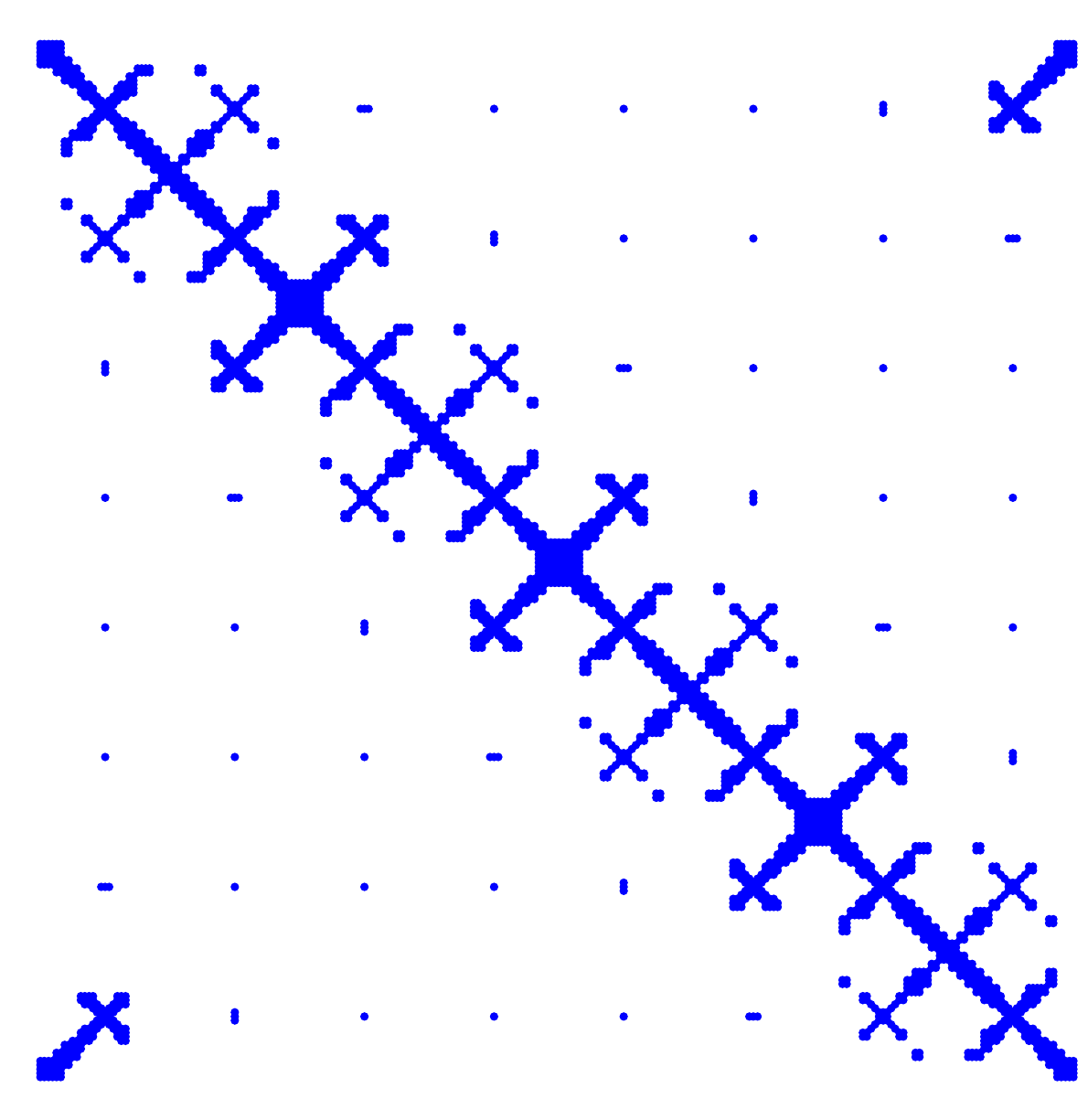}}
\fbox{\includegraphics[scale = 0.7,trim = 5.0cm 8.0cm 4.5cm 7.5cm, clip = true, scale = 0.35]{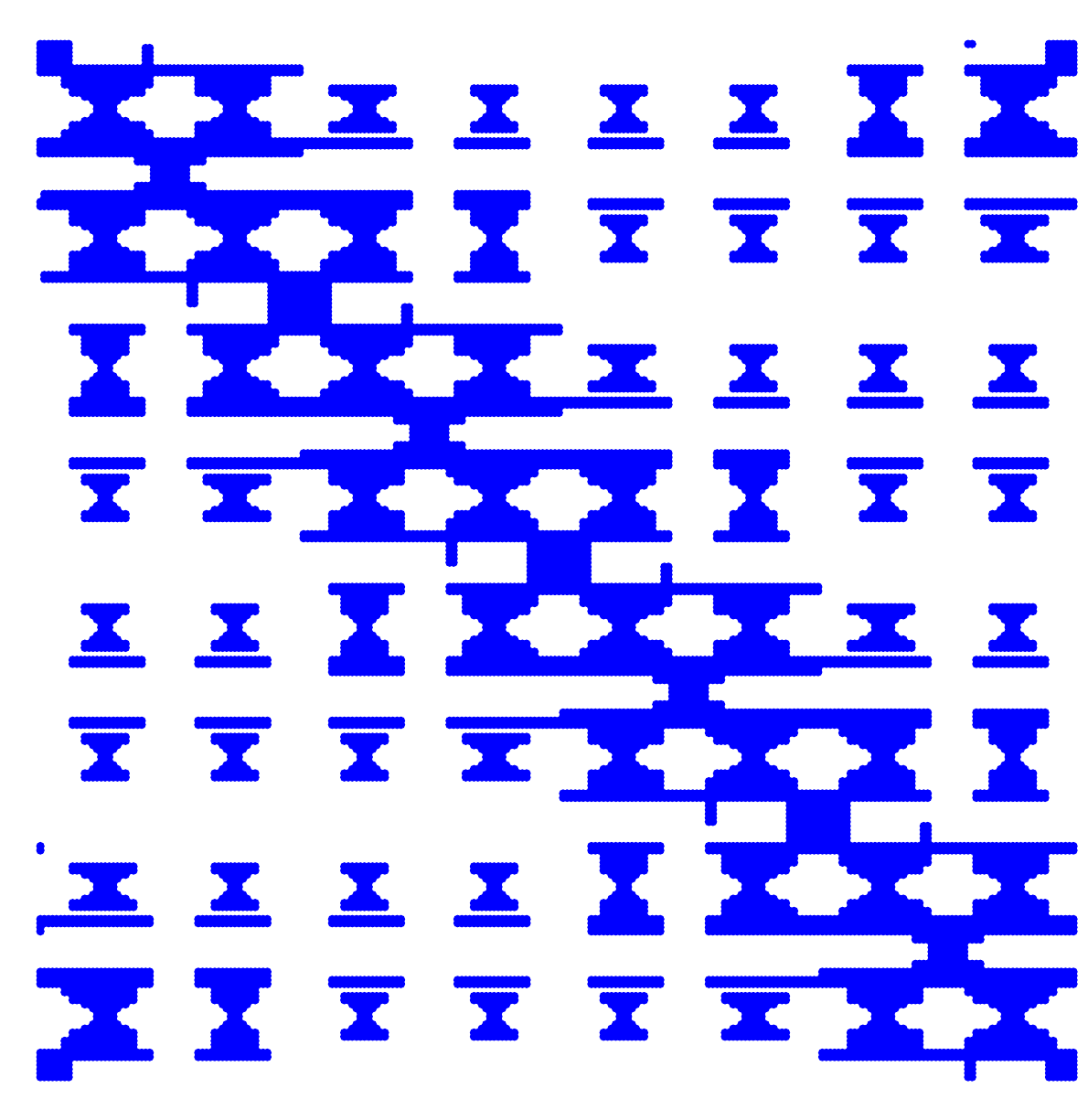}}
\caption{\label{f:sparsity} Sparsity plots of $D_{\mathrm{near}}$ for
four different preconditioners.  For these plots, we used the
eight-lobed flower from Figure~\ref{f:geom}.  The left plot corresponds
to the banded approximation $P_{B}$, the second plot corresponds to the
block diagonal approximation $P_{D}$, the next plot corresponds to the
$U$-list approximation $P_{0}$, and the final plot corresponds to the
$V_{1}$-list approximation $P_{1}$.  The small protrusions from the
diagonal of the block diagonal approximation correspond to the eight
points where the lobes of the flower approach one another.  In the
third plot, we can see that the sparsity pattern of $D_{0}$ becomes
quite complex.  A generic fast solver for such matrices is something we
are currently investigating.  The sparsity pattern of the final plot is
far too complex to construct exact inverses.}
\end{figure}
Another iterative method for solving~\eqref{e:fredholm:discrete} is
formed by first splitting $D$ into a near field and a far field.
Single-grid preconditioners can be derived by the FMM decomposition $D
= D_{\mathrm{near}}+D_{\mathrm{far}}$, where $D_{\mathrm{near}}$ is the
double-layer potential due to the near field and $D_{\mathrm{far}}$ is
the double-layer potential due to the far field.  We can iterate on the
equation
\begin{align}
  (I + D_{\mathrm{near}}) \eeta^{\mathrm{new}} = \ff - D_{\mathrm{far}}\eeta^{\mathrm{old}}.
  \label{e:nearfar:iterative}
\end{align}
This scheme can be either used as a smoother or we can use it as a
single-level preconditioner (with $\eeta^0=0$).  In the latter case, we
simply use 
\begin{align*}
  P=(I + D_{\mathrm{near}})^{-1},
\end{align*}
where the inverse indicates either an exact or an approximate
factorization. In multigrid this preconditioner can be combined
with~\eqref{e:nearfar:iterative} to define an inexact smoothing step:
\begin{align*}
  \eeta^{\mathrm{new}} = P(\ff - D_{\mathrm{far}}\eeta^{\mathrm{old}}).
\end{align*}
All the preconditioners we define below can be used in a similar manner
as approximate smoothers.

We describe several of the options for approximating
$(I+D_{\mathrm{near}})^{-1}$.
\begin{itemize}

 \item $P_{B}$, {\bf Banded approximation}: If $\Gamma$ is not too
 extreme (left two figures of Figure~\ref{f:geom}), then we can say
 that the first $s$ points to the right and left of each point are in
 the near field.  The sparsity pattern of $D_{\mathrm{near}}$ is in the
 left plot of Figure~\ref{f:sparsity}.  $P_{B}$ is very easy to compute
 and parallelize. \\

 \item $P_{D}$, {\bf Block diagonal approximation}: Alternatively, we
 can construct a quadtree and define $D_{\mathrm{near}}$ to be the
 double-layer potential due to the set of points in the same leaf $b$
 as the target point.  The sparsity pattern of $D_{\mathrm{near}}$ is
 illustrated in the second plot of Figure~\ref{f:sparsity} (this matrix
 is block diagonal after a permutation of the rows).  Since $P_{D}$ is
 block diagonal, it is very easy to compute and parallelize. \\

 \item $P_0$, {\bf $\boldsymbol{U}$-list approximation}: We can include
 more terms from the FMM in $D_{\mathrm{near}}$.  We remind the reader
 that the operators $D_{0},D_{1},D_{2},\ldots$ are the double-layer
 potential due to points in $U,V_{1},V_{2},\ldots$, respectively. The
 simplest case is to use 
 \begin{align*} 
   P_{0} \approx   A_0^{-1} = (I+D_0)^{-1},
 \end{align*} 
 where $D_{0}$ corresponds to the direct interactions in FMM. In the
 third plot of Figure~\ref{f:sparsity}, we illustrate the sparsity of
 $I + D_{0}$.  $P_{0}$ is harder to compute and parallelize. It
 requires fast sparse direct solvers. \\
 
 \item $P_\ell$, {\bf $\boldsymbol{V_{\ell}}$-list approximation}: We
 can include increasingly more and more interactions by using the
 $V_{\ell}$-list interactions. For example $P_{1} \approx (A_{0} +
 D_{1})^{-1}$.  The difficulty of constructing this preconditioner is
 even higher since $D_{1}$ is never constructed but rather is
 approximated in the FMM scheme.  The sparsity pattern of $A_{0} +
 D_{1}$ is in the right plot of Figure~\ref{f:sparsity}. \\

 \item $P_{S,\ell}$, {\bf FMMSCHUR of $\boldsymbol{P_\ell}$}: This
 preconditioner approximates $P_\ell$ using low-rank decompositions and
 the SMW formula. It is described in detail in the next section. 

\end{itemize}

These preconditioners can be used as smoothers in a classical V-cycle
multigrid scheme in which the restriction and prolongation operations
are done spectrally.  We give more details in section~\ref{s:results}.
More complicated smoothers could be built by using these
preconditioners in a splitting-based Picard or Chebyshev iterative
scheme.  However, the overall cost increases too much, and Chebyshev
smoothers require knowledge of the spectrum of $D$.  

As we mentioned in the introduction, a lot of effort has been put into
finding good approximations for $P_{0}$. However, as we will see, even
if we use an exact factorization of $(I+D_{0})$ for $P_{0}$, we still
require a lot of GMRES iterations.  What we have left out is
significant. Adding additional levels of interactions to build $P_{1}$
and $P_{2}$ cannot be done using existing technologies.  Next we
describe a new scheme, which we refer to as FMMSCHUR, to build
approximations that include additional interactions. Our scheme
requires the existence of an efficient approximation for $P_{0}$.  The
new preconditioner is termed $P_{S,L}$ and it is a single-grid
preconditioner.

\subsection{Recursive low-rank approximation preconditioner}
\label{s:schur-low-rank}
Recall that we have to invert $A=I+D$, $A \in \RR^{N\times N}$.  Also
recall the SMW formula
\begin{align*}
  (B+UV)^{-1} = (I - B^{-1}U S^{-1}V) B^{-1},\quad \mbox{where} \quad
  S=I+VB^{-1}U.
\end{align*}

First we will consider an approximation for $P_{1}$.  Define
$A_{1}=A_{0} + D_{1}$ and assume that $D_{1}$ is factorized as
$D_{1}=L_{1}M_{1}^{T}$, where $L_{1},M_{1} \in \RR^{N\times m_{1}}$,
where $m_{1}<N$.  In general, $m_{1}$ is too large to construct an
efficient preconditioner.  However, we will see in
section~\ref{s:results} that taking $m_{1}=5$ results in a good
preconditioner of $A_{1}$, even though it is a poor approximation of
$P_{1}$.  Using the SMW formula we can directly write this approximate
inverse of $A_{1}$
\begin{align*}
  A_{1}^{-1} \approx (I - A_0^{-1}L_1 S_1^{-1} M_1^T)A_0^{-1},
\end{align*}
where $S_{1} = I + M_{1}^{T} A_{0}^{-1} L_{1}$ with $S_{1} \in\RR^{m_{1} \times
m_{1}}$.

This preconditioner requires applying $S_{1}^{-1}$ and $A_{0}^{-1}$.
Since $A_{0}$ is sparse, we compute $A_{0}^{-1}$ in linear time with
the incomplete LU decomposition.  As is for the case for $D_{1} =
L_{1}M_{1}^{T}$, a rank $m_{1}=5$ approximation of the matrix
$M_{1}^{T} P_{0} L_{1} \approx UV$ results in a good preconditioner of
$S_{1}$
\begin{align*}
  \tilde{S}^{-1}_{1} = I - U (I + VU)^{-1} V.
\end{align*}
Combining all the pieces, we have constructed the preconditioner for $A_{1}$
\begin{align}
  P_{S,1} = (I - P_{0} L_{1} \tilde{S}^{-1}_{1} M_{1}^{T})P_{0}, 
  \label{e:schur:preco}
\end{align}
where $P_{0}$ and $ \tilde{S}^{-1}_{1}$ are preconditioners for $A_{0}$
and $S_{1}$, respectively.  Note that applying this preconditioner
requires applying $P_{0}$ twice.

The cost of constructing $ \tilde{S}^{-1}_{1}$ depends on how we
compute the low-rank approximation of $M_{1}^{T} P_{0} L_{1}$. Using
techniques discussed in~\cite{halko-martinsson-tropp11}, the low-rank
vectors can be approximated in work that is proportional to the work of
applying $P_{0}$.  These techniques for forming the require low-rank
approximations are left as future work, and in this work, we simply
use the truncated SVD.

Now consider the case of $A_{2}=A_{1} + D_{2}$.  The scheme can be
applied recursively in which the approximate factorization of $A_{1}$
described above can be used to build a preconditioner for $A_{2}$
\begin{align*}
  P_{S,2} = (I - P_{S,1} L_{2}  \tilde{S}^{-1}_{2} M_{2}^{T})P_{S,1}, 
\end{align*}
where $ \tilde{S}^{-1}_{2} \approx (I + M_{2}^{T}P_{S,1}L_{2})^{-1}$;
$\tilde{S}^{-1}_{2}$ is built by computing a low-rank approximation of
$M_{2}^{T}P_{S,1}L_{2}$.  Applying this preconditioner requires four
applications of $P_{0}$.  This scheme can be used recursively to add
more levels. The main observation is that in our scheme $D_{\ell}$ has
a smaller norm and stronger compactness than $D$ since it only involves
the far field and no kernel singularities. So the Schur complements
$S_{\ell} = I + R_{\ell}$ will be very well-conditioned matrices, with
$R_{\ell}$ admitting low-rank approximations that can be used to
construct good single-grid preconditioners.

Unfortunately, if we use recursion, the cost of applying $P_{S,L}$ is
$2^{L}$ times the cost of applying $P_0$.  Since $L=\bigO(\log_4(N/s))$
in 2D, the cost of applying $P_{S,L}$ is $\bigO(N^{3/2})$. Therefore we
need to consider either alternative algorithms that do not apply
$P_{S,L}$ recursively or use small $L$, say $L=1$. In this paper we
only consider the latter case.  If we use small $L$, the preconditioner
can be further improved by approximating the effect of all other levels
by a global low-rank approximation.  That is, we write
\begin{align*}  
  A = A_0 + \sum_{i=1}^L D_\ell + Z.
\end{align*}
We build $P_{S,L}$ to approximate the inverse of $A_0 +
\sum_{i=1}^LD_\ell$ and we use a low-rank approximation of $Z$, which
we combine with $P_{S,L}$ using the SMW formula. The number of vectors
used in approximating $Z$ is selected so that the cost of applying
$P_{S,L}$ is $\bigO(N\log N)$; we choose it to be $\bigO(\log N)$.  In
the next section we provide experimental evidence that using $L=1$ with
a $\bigO(\log N)$-rank approximation of $Z$ works quite well. 

To see the connection to the FMM, recall the derivation of the SMW
formula is based on introducing auxiliary unknowns and using a Schur
complement decomposition. Indeed if we want to solve $(A+D_1)\eeta =
f$, we set $\eeta_{1} = M_{1}^{T}\eeta$ and solve a Schur complement
for $\eeta_{1}$ instead of the original system.  In the FMM context,
$\eeta_{1}$ corresponds to the multipole moments $M_{1}^{T} \eeta$.
For this reason we use the subscript ``S'' to indicate that it is based
on Schur complement decomposition.

As mentioned in the introduction, this preconditioner is different than
the fast algorithms for integral equations. Roughly speaking, those are
based on non-overlapping Schur decompositions of the domain combined
with specially constructed low-rank approximations of the interface and
tailored to boundary integral equations. Our scheme is less powerful
but can be implemented relatively easily on top of a standard FMM
evaluation.

\section{Results}
\label{s:results}

We discuss the behavior of the different preconditioners in the context
of single-grid and multigrid schemes. We conduct a series of
experiments, which we now summarize.

\begin{itemize}
\item {\bf Single-grid and two-level multigrid}
  (Table~\ref{t:twoLevel}): Here we compare the effectiveness of the
  preconditioners in a single-grid and a two-grid scheme applied to
  several geometries. In the single-grid case the new FMM-based
  preconditioners are more effective.  We see that if the coarse grid
  is sufficiently resolved, then the Picard and $P_0$ preconditioners
  work best.

\item {\bf Geometric versus projection coarsening}
  (Tables~\ref{t:curly:gmres}--\ref{t:8flower:gmres}): Here we compare
  the cost of the two methods for creating the coarse grid operators.
  The projection-based coarse grid operator uses fewer GMRES steps,
  but it requires too many matvecs.  We again see the effect of using
  unresolved coarse grids.  Also, the simpler preconditioners $P_{B}$
  and $P_{D}$ do not reduce the GMRES iterations significantly. 

\item {\bf Smoothing properties and overall multigrid performance}
  (Table~\ref{t:coarse512} and Figure~\ref{f:errSmoother}): Here we
  further analyze the different preconditioners and consider two
  different coarse grids applied to the four-lobed flower.  We see
  that the Picard preconditioner does the best job of eliminating the
  high frequencies at multiple resolutions.  However, at low
  resolutions, the Picard smoother is divergent and the resulting
  preconditioner is ineffective.  The FMM-based smoothers are
  convergent and the resulting preconditioner reduces the total number
  of GMRES iterations.  We also demonstrate that post-smoothing does
  not sufficiently reduce the number of GMRES iterations to justify
  its use.

\item {\bf Single-grid preconditioners for unresolved geometries}
  (Table~\ref{t:psc-preco}): We assess the single-grid preconditioners
  for two problems.  A strategy that couples a Picard two-grid
  preconditioner with a single-grid FMM-based preconditioner is
  outlined.  We apply several preconditioners to the 8-lobed and
  24-lobed flowers.  We see that FMMSCHUR preconditioner $P_{S,1}$
  works best which indicates that it is important to include part the
  far field in the preconditioner, even if it is a very crude low-rank
  approximation.
\end{itemize}

Other details of our results include:
\begin{itemize}

\item We consider solving~\eqref{e:fredholm:discrete} for the four
geometries in Figure~\ref{f:geom}.  These geometries are selected so
that the number of points required to represent the geometry range from
a few hundred to several thousands of points.

\item Our implementation is in MATLAB.  Some of the operators and
factorizations are not computed by fast algorithms, so the actual
wall-clock times are not informative.  For example, all low-rank
approximations are computed by truncating MATLAB's exact SVD to the
desired rank. All factorizations of the preconditioners are exact
unless otherwise noted. When we use incomplete factorizations, we use
MATLAB's incomplete LU factorization (ILU) with a drop tolerance of
$10^{-3}$.  To compare the different methodologies, we estimate the
cost of the preconditioner in terms of the cost to apply $D$ using the
FMM.

\item In all of our runs, we use GMRES without restarts preconditioned
with the single-grid schemes we discussed in section~\ref{s:preco} and
V-cycle multigrid.  GMRES is terminated when the relative tolerance
dropped below $10^{-12}$. In practice, there is no reason to have an
algebraic error that is smaller than the truncation error.  But for
simplicity, we keep the GMRES tolerance fixed in all of our tests.

\item In GMRES, we always use the exact double-layer potential; that
is, we do not use any FMM approximation.  FMM is used only to build the
preconditioners.  In all runs, unless otherwise stated, we only use
four moments for the upward, downward, and translation operators.  That
is, the far field is approximated with four degrees of freedom per
box.  In this way, the far field is approximated as crudely as possible
within the FMM context.

\item For the preconditioners $P_{\ell}$, we have results with $s=4$,
$s=10$, and $s=50$ points per leaf.  Although we refer to $P_{1}$ as
including the $V_{1}$-list interactions, it also includes other
interactions since the trees are adaptive (Figure~\ref{f:quadtrees}
and~\ref{f:mgLayers}).  In particular, $P_{1}$ also includes the
interactions due to the $W$-list and $X$-list, where these lists are
defined in~\cite{car:gre:rok}.

\item For the preconditioner $P_{S,1}$, we use a rank five
approximation for $D_{1} = L_{1}M_{1}^{T}$ and for the Schur complement
$M_{1}^TP_{0}L_{1}$, and a rank $5\log(N/s)$ approximation for $Z =
A-A_{0}$.  We have conducted experiments for $P_{S,1}$ where we set
$P_{0}$ with an ILU factorization of $A_0$.  We discuss this further in
section~\ref{s:singleGrid}.

\item Additional parameters include the number of pre- and
post-smoothing steps and the number of levels.  Results are presented
for single-grid, two-grid, and multigrid preconditioners.

\end{itemize}

\paragraph{Estimation of the preconditioner cost.}
Before we present the results in detail, we discuss how we estimate the
cost of each preconditioner.  We report the number of scaled matvecs
required by our unpreconditioned and different preconditioned runs.
These are scaled assuming that the matrix $D$ can be applied in linear
time.  For $P_{B},P_{D},P_{0}$, and $P_{1}$, we assume that their
application costs the same as a matvec with $D$.  Since $P_{B}$ and
$P_{D}$ have a predictable structure, they can be applied efficiently
with standard banded and block diagonal solvers.  However, our
assumption for $P_{0}$ and $P_{1}$ is quite strong since they require
exact factorizations of a sparse but unstructured matrix.  Fortunately,
an ILU decomposition of $P_{0}$, which can be computed in linear time,
can be used without any major effect on the quality of the
preconditioner.  Concerning $P_{1}$, its exact factorization is
impossible and this motivates the single-grid preconditioner
$P_{S,1}$.  The cost of applying $P_{S,1}$ is $\bigO(N \log N)$ since
we need to approximate $Z$. For the cost of multigrid see
section~\ref{s:geo-vs-proj}.

\subsection{Single-grid and two-level multigrid, Table~\ref{t:twoLevel}}
\begin{table}
\centering
\begin{minipage}[t]{0.12\textwidth}
\vspace{40pt}
\centering
\ifInputs
  \scalebox{0.4}{\input{curly.tikz}}
\fi
\end{minipage}
\begin{minipage}[t]{0.65\textwidth}
\vspace{0pt}
\centering
\vspace{5pt}
$N = 2048$ \\
\vspace{5pt}
\begin{tabular}{cccc}
\hline
Preconditioner & Single-Grid & $N_{\min} = 128$ & $N_{\min} = 16$ \\
\hline
Picard & 16 & 3 & 13 \\
$P_{D}$ & 16 & 8 & 13 \\
$P_{0}$ & 16 & 8 & 13 \\
$P_{1}$-Exact & 15 & 8 & 13 \\
$P_{1}$ & 15 & 8 & 13 \\
$P_{2}$ & 14 & 7 & 12 \\
$P_{S,1}$ & 7 & - & -\\
\hline
\end{tabular}
\end{minipage}
\begin{minipage}[t]{0.21\textwidth}
\vspace{30pt}
\centering
\ifInputs
  \scalebox{0.75}{\begin{tikzpicture}[scale = 0.5]

\begin{axis}[
  xmin = 0,
  xmax = 17,
  ymin = 1e-14,
  ymax = 60,
  xlabel = {Iteration},
  ylabel = {Relative Residual},
  ymode = log,
  xtick = {0,5,10,15},
  xticklabels = {\Large{0},\Large{5},\Large{10},\Large{15}},
  ytick = {1e-12,1e-6,1e0},
  yticklabels = {\Large{$10^{-12}$},\Large{$10^{-6}$},\Large{$10^{0}$}},
  scaled x ticks = false,
  scaled y ticks = false,
  label style = {font=\Large},
  ylabel style = {yshift = 10pt},
  grid,
  legend cell align = left,
  legend entries = {Unpre,Picard,$P_{D}$,$P_{0}$,$P_{1}$},
  legend pos = north east,
  legend style = {draw=none}
]

\addplot[mark=none,blue,mark size = 1pt,line width = 1.5pt] table[col sep=comma]{ 
1,1.0000e+00
2,1.4392e-01
3,7.3199e-02
4,2.1202e-02
5,5.2634e-03
6,1.2735e-03
7,2.0351e-04
8,3.6000e-05
9,5.7257e-06
10,1.0737e-06
11,1.8982e-07
12,2.8160e-08
13,5.3924e-09
14,5.0874e-10
15,7.1254e-11
16,1.1265e-11
17,9.6012e-13
}; 

\addplot[mark=none,red,mark size = 1pt,line width = 1.5pt] table[col sep=comma] {
1,2.2440e+00
2,5.1622e-05
3,6.4040e-09
4,2.9922e-13
};

\addplot[mark=none,green,mark size = 1pt,line width = 1.5pt] table[col sep=comma] {
1,2.2440e+00
2,1.1731e-02
3,3.0327e-04
4,8.1306e-06
5,2.2704e-07
6,7.8820e-09
7,2.7419e-10
8,6.0059e-12
9,1.1429e-13
};

\addplot[mark=none,black,mark size = 1pt,line width = 1.5pt] table[col sep=comma] {
1,2.2444e+00
2,4.2686e-02
3,1.7175e-03
4,4.7275e-05
5,8.9444e-07
6,1.9954e-08
7,3.8876e-10
8,8.5878e-12
9,2.0645e-13
};

\addplot[mark=none,cyan,mark size = 1pt,line width = 1.5pt] table[col sep=comma] {
1,2.2450e+00
2,6.8606e-02
3,1.3130e-03
4,1.6326e-05
5,3.7779e-07
6,8.4074e-09
7,2.0297e-10
8,3.0082e-12
9,4.7610e-14
};
\end{axis}

\end{tikzpicture}}
\fi
\end{minipage}

\vspace{10pt}
\begin{minipage}[t]{0.12\textwidth}
\vspace{40pt}
\centering
\ifInputs
  \scalebox{0.4}{\begin{tikzpicture}

\begin{axis}[
  width=1.2in,height=2in,
  scale only axis,
  xmin = -1.0,xmax = 1.6,
  ymin = -1.5,ymax = 1.5,
  hide axis,
  axis equal = true
  ]
\addplot[color=red,line width = 1.0pt,solid,domain=0:360,samples=360*2]({(1+0.5*cos(3*\x)+0.05*cos(30*\x))*cos(\x)},{(1+0.5*cos(3*\x)+0.05*cos(30*\x))*sin(\x)});
\end{axis}

\end{tikzpicture}}
\fi
\end{minipage}
\begin{minipage}[t]{0.65\textwidth}
\vspace{0pt}
\centering
$N = 2048$ \\
\vspace{5pt} 
\begin{tabular}{cccc}
\hline
Preconditioner & Single-Grid & $N_{\min} = 128$ & $N_{\min} = 16$ \\
\hline
Picard & 26 & 13 & 33 \\
$P_{D}$ & 27 & 15 & 31 \\
$P_{0}$ & 25 & 15 & 25 \\
$P_{1}$-Exact & 20 & 11 & 17 \\
$P_{1}$ & 20 & 11 & 17 \\
$P_{2}$ & 17 & 10 & 14 \\
$P_{S,1}$ & 9 & - & -\\
\hline
\end{tabular}
\end{minipage}
\begin{minipage}[t]{0.21\textwidth}
\vspace{30pt}
\centering
\ifInputs
  \scalebox{0.75}{\begin{tikzpicture}[scale = 0.5]

\begin{axis}[
  xmin = 0,
  xmax = 32,
  ymin = 1e-14,
  ymax = 60,
  xlabel = {Iteration},
  ylabel = {Relative Residual},
  ymode = log,
  xtick = {0,10,20,30},
  xticklabels = {\Large{0},\Large{5},\Large{10},\Large{20}},
  ytick = {1e-12,1e-6,1e0},
  yticklabels = {\Large{$10^{-12}$},\Large{$10^{-6}$},\Large{$10^{0}$}},
  scaled x ticks = false,
  scaled y ticks = false,
  label style = {font=\Large},
  ylabel style = {yshift = 10pt},
  grid,
  legend cell align = left,
  legend entries = {Unpre,Picard,$P_{D}$,$P_{0}$,$P_{1}$},
  legend pos = north east,
  legend style = {draw=none}
]

\addplot[mark=none,blue,mark size = 1pt,line width = 1.5pt] table[col sep=comma] {
1,1.0000e+00
2,2.4874e-01
3,1.0379e-01
4,4.2172e-02
5,2.2328e-02
6,6.8207e-03
7,2.8317e-03
8,1.3335e-03
9,8.2225e-04
10,4.2167e-04
11,1.3300e-04
12,4.1945e-05
13,7.9343e-06
14,2.5818e-06
15,6.4557e-07
16,1.8357e-07
17,8.3003e-08
18,3.7777e-08
19,1.6127e-08
20,9.1688e-09
21,2.5970e-09
22,5.7306e-10
23,1.4860e-10
24,3.6124e-11
25,8.4015e-12
26,2.3960e-12
27,3.9104e-13
};

\addplot[mark=none,red,mark size = 1pt,line width = 1.5pt] table[col sep=comma] 
{
1,1.8300e+00
2,1.6176e-01
3,2.1631e-02
4,1.8716e-03
5,2.6747e-04
6,3.7343e-05
7,4.4022e-06
8,8.3427e-07
9,1.3787e-07
10,1.4591e-08
11,1.4173e-09
12,1.3927e-10
13,1.4953e-11
14,1.0769e-12
};

\addplot[mark=none,green,mark size = 1pt,line width = 1.5pt] table[col sep=comma] 
{
1,1.8352e+00
2,1.9063e-01
3,3.1721e-02
4,4.3595e-03
5,9.8320e-04
6,1.4742e-04
7,2.6714e-05
8,4.9337e-06
9,8.2231e-07
10,1.0973e-07
11,1.5898e-08
12,2.6228e-09
13,3.4354e-10
14,3.8869e-11
15,5.4027e-12
16,6.4473e-13
};

\addplot[mark=none,black,mark size = 1pt,line width = 1.5pt] table[col sep=comma] 
{
1,1.8483e+00
2,2.3271e-01
3,3.3497e-02
4,4.4976e-03
5,5.9381e-04
6,9.8597e-05
7,1.4647e-05
8,1.8734e-06
9,2.5199e-07
10,2.5639e-08
11,3.7414e-09
12,6.6256e-10
13,1.2401e-10
14,1.7778e-11
15,2.9515e-12
16,4.2403e-13
};

\addplot[mark=none,cyan,mark size = 1pt,line width = 1.5pt] table[col sep=comma] 
{
1,1.8334e+00
2,2.4259e-01
3,2.5409e-02
4,1.6271e-03
5,1.5704e-04
6,1.0408e-05
7,6.3331e-07
8,3.1292e-08
9,1.9083e-09
10,1.0064e-10
11,5.2556e-12
12,2.8835e-13
};

\end{axis}

\end{tikzpicture}}
\fi
\end{minipage}

\vspace{10pt}
\begin{minipage}[t]{0.12\textwidth}
\vspace{40pt}
\centering
\ifInputs
  \scalebox{0.4}{\begin{tikzpicture}

\begin{axis}[
  width=2in, height=2in,
  scale only axis,
  xmin=-2, xmax=2,
  ymin=-2, ymax=2,
  hide axis
  ]
\addplot[color=red,line width = 1.0pt,solid,domain=0:360,samples=360*2]({(1+0.98*cos(4*\x))*cos(\x)},{(1+0.98*cos(4*\x))*sin(\x)});
\end{axis}

\end{tikzpicture}}
\fi
\end{minipage}
\begin{minipage}[t]{0.65\textwidth}
\vspace{0pt}
\centering
$N = 2048$ \\
\vspace{5pt}
\begin{tabular}{cccc}
\hline
Preconditioner & Single-Grid & $N_{\min} = 128$ & $N_{\min} = 16$ \\
\hline
Picard & 42 & 55 & 55 \\
$P_{D}$ & 44 & 46 & 55 \\
$P_{0}$ & 30 & 25 & 31 \\
$P_{1}$-Exact & 23 & 15 & 23 \\
$P_{1}$ & 23 & 15 & 23 \\
$P_{2}$ & 18 & 13 & 19 \\
$P_{S,1}$ & 10 & - & -\\
\hline
\end{tabular}
\end{minipage}
\begin{minipage}[t]{0.21\textwidth}
\vspace{30pt}
\centering
\ifInputs
  \scalebox{0.75}{\begin{tikzpicture}[scale = 0.5]

\begin{axis}[
  xmin = 0,
  xmax = 55,
  ymin = 1e-14,
  ymax = 60,
  xlabel = {Iteration},
  ylabel = {Relative Residual},
  ymode = log,
  xtick = {0,25,50},
  xticklabels = {\Large{0},\Large{25},\Large{50}},
  ytick = {1e-12,1e-6,1e0},
  yticklabels = {\Large{$10^{-12}$},\Large{$10^{-6}$},\Large{$10^{0}$}},
  scaled x ticks = false,
  scaled y ticks = false,
  label style = {font=\Large},
  ylabel style = {yshift = 10pt},
  grid,
  legend cell align = left,
  legend entries = {Unpre,Picard,$P_{D}$,$P_{0}$,$P_{1}$},
  legend pos = north east,
  legend style = {draw=none}
]

\addplot[mark=none,blue,mark size = 1pt,line width = 1.5pt] table[col sep=comma]
{
1,1.0000e+00
2,3.0301e-01
3,8.8663e-02
4,4.4853e-02
5,3.9284e-02
6,2.3072e-02
7,9.8487e-03
8,7.2061e-03
9,6.3034e-03
10,5.1127e-03
11,4.0781e-03
12,3.6853e-03
13,2.7593e-03
14,2.0605e-03
15,1.1834e-03
16,6.2715e-04
17,3.3116e-04
18,1.6999e-04
19,1.1619e-04
20,8.0046e-05
21,5.8625e-05
22,3.6583e-05
23,2.1395e-05
24,1.1618e-05
25,6.4305e-06
26,3.0433e-06
27,1.3234e-06
28,5.6248e-07
29,2.4340e-07
30,1.1518e-07
31,4.9420e-08
32,2.2917e-08
33,9.8886e-09
34,4.4733e-09
35,1.6804e-09
36,5.7976e-10
37,2.0840e-10
38,7.9111e-11
39,3.3307e-11
40,1.1103e-11
41,3.5980e-12
42,1.0423e-12
43,3.4100e-13
};

\addplot[mark=none,red,mark size = 1pt,line width = 1.5pt] table[col sep=comma]
{
1,1.9188e+00
2,1.0995e+00
3,2.6005e-01
4,1.5623e-01
5,1.4744e-01
6,3.8291e-02
7,2.7322e-02
8,2.5415e-02
9,2.1235e-02
10,1.5968e-02
11,9.4813e-03
12,5.8569e-03
13,4.2789e-03
14,3.4876e-03
15,3.1643e-03
16,2.9891e-03
17,2.7389e-03
18,2.4858e-03
19,2.3656e-03
20,2.2508e-03
21,2.1723e-03
22,1.9705e-03
23,1.7090e-03
24,1.3573e-03
25,8.8635e-04
26,5.4204e-04
27,2.5977e-04
28,1.3657e-04
29,6.1390e-05
30,3.4538e-05
31,2.6618e-05
32,2.2593e-05
33,1.4643e-05
34,9.2043e-06
35,9.1932e-06
36,5.8357e-06
37,2.8346e-06
38,1.5404e-06
39,9.4542e-07
40,4.9145e-07
41,3.0717e-07
42,3.0468e-07
43,2.0412e-07
44,9.4564e-08
45,5.3384e-08
46,3.3056e-08
47,1.1525e-08
48,2.6672e-09
49,9.8745e-10
50,4.2502e-10
51,1.3482e-10
52,4.8705e-11
53,2.1015e-11
54,7.3473e-12
55,2.9969e-12
56,9.9279e-13
};

\addplot[mark=none,green,mark size = 1pt,line width = 1.5pt] table[col sep=comma]
{
1,1.9236e+00
2,8.9195e-01
3,4.4174e-01
4,1.4849e-01
5,1.4667e-01
6,9.0906e-02
7,6.7754e-02
8,3.3037e-02
9,2.1609e-02
10,1.6733e-02
11,1.0839e-02
12,9.1695e-03
13,7.8716e-03
14,6.3344e-03
15,3.9738e-03
16,2.0151e-03
17,1.4915e-03
18,1.4250e-03
19,7.2896e-04
20,2.9968e-04
21,2.7888e-04
22,2.0789e-04
23,8.8627e-05
24,5.8033e-05
25,2.4802e-05
26,9.1886e-06
27,3.5557e-06
28,1.1598e-06
29,4.6247e-07
30,3.1924e-07
31,2.4208e-07
32,2.2462e-07
33,1.9348e-07
34,6.8477e-08
35,1.9626e-08
36,1.1574e-08
37,9.3570e-09
38,8.5533e-09
39,5.6206e-09
40,1.9029e-09
41,6.8674e-10
42,3.9239e-10
43,2.3045e-10
44,5.9788e-11
45,1.4397e-11
46,3.4048e-12
47,1.1377e-12
};

\addplot[mark=none,black,mark size = 1pt,line width = 1.5pt] table[col sep=comma]
{
1,2.0825e+00
2,7.6275e-01
3,1.5511e-01
4,6.3016e-02
5,1.9940e-02
6,7.1449e-03
7,4.2650e-03
8,3.3801e-03
9,8.3692e-04
10,1.2201e-04
11,1.8494e-05
12,3.6404e-06
13,1.8511e-06
14,1.6141e-06
15,8.6423e-07
16,3.5289e-07
17,1.0654e-07
18,2.1364e-08
19,4.0564e-09
20,2.6621e-09
21,2.1741e-09
22,1.0990e-09
23,1.3529e-10
24,2.4141e-11
25,7.1105e-12
26,1.0148e-12
};

\addplot[mark=none,cyan,mark size = 1pt,line width = 1.5pt] table[col sep=comma]
{
1,2.5294e+00
2,1.3783e+00
3,2.6911e-01
4,1.2357e-01
5,3.9246e-02
6,3.0941e-03
7,2.4734e-04
8,1.0604e-05
9,1.8981e-06
10,1.1130e-06
11,1.5819e-07
12,4.3887e-08
13,7.8241e-09
14,4.6198e-10
15,4.6799e-11
16,1.4085e-12
};

\end{axis}

\end{tikzpicture}}
\fi
\end{minipage}

\vspace{10pt}
\begin{minipage}[t]{0.12\textwidth}
\vspace{40pt}
\centering
\ifInputs
  \scalebox{0.4}{\begin{tikzpicture}

\begin{axis}[
  width=2in, height=2in,
  scale only axis,
  xmin=-2, xmax=2,
  ymin=-2, ymax=2,
  hide axis
  ]
\addplot[color=red,line width = 1.0pt,solid,domain=0:360,samples=360*2]({(1+0.98*cos(8*\x))*cos(\x)},{(1+0.98*cos(8*\x))*sin(\x)});
\end{axis}

\end{tikzpicture}}
\fi
\end{minipage}
\begin{minipage}[t]{0.65\textwidth}
\vspace{0pt}
\centering
$N = 4096$ \\
\vspace{5pt}
\begin{tabular}{cccc}
\hline
Preconditioner & Single-Grid & $N_{\min} = 128$ & $N_{\min} = 16$ \\
\hline
Picard & 115 & 190 & 179 \\
$P_{D}$ & 99 & 140 & 125 \\
$P_{0}$ & 55 & 70 & 60 \\
$P_{1}$-Exact & 39 & 33 & 47 \\
$P_{1}$ & 39 & 33 & 47 \\
$P_{2}$ & 27 & 26 & 41 \\
$P_{S,1}$ & 19 & - & -\\
\hline
\end{tabular}
\end{minipage}
\begin{minipage}[t]{0.21\textwidth}
\vspace{30pt}
\centering
\ifInputs
  \scalebox{0.75}{\begin{tikzpicture}[scale = 0.5]

\begin{axis}[
  xmin = 0,
  xmax = 210,
  ymin = 1e-14,
  ymax = 60,
  xlabel = {Iteration},
  ylabel = {Relative Residual},
  ymode = log,
  xtick = {0,100,200},
  xticklabels = {\Large{0},\Large{100},\Large{200}},
  ytick = {1e-12,1e-6,1e0},
  yticklabels = {\Large{$10^{-12}$},\Large{$10^{-6}$},\Large{$10^{0}$}},
  scaled x ticks = false,
  scaled y ticks = false,
  label style = {font=\Large},
  ylabel style = {yshift = 10pt},
  grid,
  legend cell align = left,
  legend entries = {Unpre,Picard,$P_{D}$,$P_{0}$,$P_{1}$},
  legend pos = north east,
  legend style = {draw=none}
]

\addplot[mark=none,blue,mark size = 1pt,line width = 1.5pt] table[col sep=comma]
{
1,1.0000e+00
2,1.5700e-01
3,2.5951e-02
4,1.8296e-02
5,1.4714e-02
6,1.2740e-02
7,1.2479e-02
8,1.0924e-02
9,8.8122e-03
10,6.3501e-03
11,4.3262e-03
12,3.2151e-03
13,2.2760e-03
14,1.7374e-03
15,1.4743e-03
16,1.3435e-03
17,1.3312e-03
18,1.3226e-03
19,1.2203e-03
20,1.1013e-03
21,9.1759e-04
22,7.9336e-04
23,6.4184e-04
24,5.2348e-04
25,4.0693e-04
26,3.2779e-04
27,2.6306e-04
28,2.3512e-04
29,2.2093e-04
30,2.1937e-04
31,2.1860e-04
32,2.1165e-04
33,1.9822e-04
34,1.8438e-04
35,1.6477e-04
36,1.4743e-04
37,1.2770e-04
38,1.1746e-04
39,1.0992e-04
40,1.0707e-04
41,1.0454e-04
42,1.0324e-04
43,1.0177e-04
44,1.0172e-04
45,1.0163e-04
46,1.0076e-04
47,9.6568e-05
48,9.1427e-05
49,8.2042e-05
50,7.6445e-05
51,7.1522e-05
52,6.6362e-05
53,6.2247e-05
54,5.4405e-05
55,4.9466e-05
56,4.1200e-05
57,3.5672e-05
58,2.8516e-05
59,2.4658e-05
60,2.1151e-05
61,2.0004e-05
62,1.8429e-05
63,1.7808e-05
64,1.6631e-05
65,1.5364e-05
66,1.3535e-05
67,1.0353e-05
68,7.4944e-06
69,5.0773e-06
70,3.5850e-06
71,2.4465e-06
72,1.7238e-06
73,1.1194e-06
74,7.6857e-07
75,4.5590e-07
76,2.7529e-07
77,1.6133e-07
78,1.0200e-07
79,7.2967e-08
80,5.4322e-08
81,4.2263e-08
82,3.2197e-08
83,2.1813e-08
84,1.3961e-08
85,7.9153e-09
86,4.8622e-09
87,3.5867e-09
88,3.0593e-09
89,3.0169e-09
90,2.9574e-09
91,2.5171e-09
92,1.8892e-09
93,1.3432e-09
94,9.2484e-10
95,6.8514e-10
96,5.6472e-10
97,5.1489e-10
98,4.9701e-10
99,4.9699e-10
100,4.8068e-10
101,4.1123e-10
102,3.0297e-10
103,1.8631e-10
104,1.1226e-10
105,6.5871e-11
106,3.9813e-11
107,2.3448e-11
108,1.4625e-11
109,9.8273e-12
110,9.0783e-12
111,9.0166e-12
112,7.5463e-12
113,4.6070e-12
114,2.5069e-12
115,1.2604e-12
116,6.9279e-13
};

\addplot[mark=none,red,mark size = 1pt,line width = 1.5pt] table[col sep=comma]
{
1,1.8338e+00
2,1.7844e+00
3,5.7970e-01
4,1.2807e-01
5,5.3880e-02
6,4.4645e-02
7,4.1676e-02
8,2.2441e-02
9,2.0940e-02
10,2.0567e-02
11,1.8772e-02
12,1.4641e-02
13,1.2346e-02
14,1.1413e-02
15,1.1052e-02
16,1.0827e-02
17,1.0778e-02
18,9.9467e-03
19,9.3953e-03
20,7.7235e-03
21,7.0775e-03
22,6.2341e-03
23,5.2623e-03
24,4.8150e-03
25,4.7370e-03
26,4.6821e-03
27,4.6490e-03
28,4.4873e-03
29,4.1697e-03
30,3.6036e-03
31,2.9752e-03
32,2.6158e-03
33,2.2463e-03
34,1.9336e-03
35,1.6421e-03
36,1.4551e-03
37,1.2814e-03
38,1.1483e-03
39,1.0498e-03
40,9.3120e-04
41,8.7330e-04
42,8.0626e-04
43,7.6142e-04
44,6.9681e-04
45,6.6295e-04
46,6.3643e-04
47,6.2913e-04
48,6.2865e-04
49,6.2826e-04
50,6.1490e-04
51,5.8777e-04
52,5.4867e-04
53,4.8900e-04
54,4.3401e-04
55,3.8811e-04
56,3.5301e-04
57,3.3245e-04
58,3.2934e-04
59,3.2934e-04
60,3.1429e-04
61,3.0034e-04
62,2.8341e-04
63,2.6336e-04
64,2.4573e-04
65,2.2063e-04
66,1.9642e-04
67,1.7314e-04
68,1.5379e-04
69,1.4006e-04
70,1.3182e-04
71,1.2455e-04
72,1.2110e-04
73,1.1800e-04
74,1.1508e-04
75,1.1312e-04
76,1.1131e-04
77,1.1103e-04
78,1.1094e-04
79,1.1089e-04
80,1.1024e-04
81,1.0958e-04
82,1.0795e-04
83,1.0614e-04
84,1.0507e-04
85,1.0314e-04
86,9.9555e-05
87,9.7951e-05
88,9.0512e-05
89,8.4747e-05
90,7.7578e-05
91,7.2791e-05
92,6.9572e-05
93,6.5524e-05
94,6.2260e-05
95,6.0307e-05
96,5.7340e-05
97,5.6625e-05
98,5.5434e-05
99,5.5202e-05
100,5.5184e-05
101,5.5181e-05
102,5.3754e-05
103,5.2609e-05
104,5.1650e-05
105,5.0896e-05
106,5.0791e-05
107,5.0781e-05
108,4.9251e-05
109,4.7397e-05
110,4.5469e-05
111,4.0985e-05
112,3.7945e-05
113,3.4651e-05
114,3.3062e-05
115,3.1877e-05
116,3.0932e-05
117,3.0285e-05
118,2.9468e-05
119,2.8824e-05
120,2.7973e-05
121,2.7014e-05
122,2.5295e-05
123,2.3654e-05
124,2.1770e-05
125,1.9772e-05
126,1.7190e-05
127,1.5024e-05
128,1.3082e-05
129,1.0858e-05
130,8.5545e-06
131,6.6243e-06
132,5.2524e-06
133,4.4995e-06
134,4.3454e-06
135,4.3316e-06
136,4.0330e-06
137,3.2646e-06
138,2.4218e-06
139,1.5477e-06
140,1.0165e-06
141,7.1835e-07
142,5.0822e-07
143,3.7215e-07
144,2.9111e-07
145,2.2150e-07
146,1.5432e-07
147,1.1677e-07
148,9.8012e-08
149,9.3582e-08
150,9.3204e-08
151,8.4690e-08
152,6.5981e-08
153,4.5033e-08
154,2.9657e-08
155,1.9506e-08
156,1.4846e-08
157,1.2152e-08
158,1.0942e-08
159,1.0259e-08
160,9.4826e-09
161,8.5551e-09
162,7.5025e-09
163,5.7678e-09
164,4.1051e-09
165,2.6626e-09
166,1.5272e-09
167,9.9860e-10
168,6.9288e-10
169,4.9148e-10
170,3.9189e-10
171,3.0727e-10
172,2.2444e-10
173,1.6532e-10
174,1.1517e-10
175,8.0524e-11
176,5.7336e-11
177,4.0945e-11
178,3.1600e-11
179,2.6485e-11
180,2.5009e-11
181,2.4975e-11
182,2.4625e-11
183,2.2660e-11
184,1.7532e-11
185,1.2155e-11
186,8.4279e-12
187,5.8499e-12
188,4.2693e-12
189,3.0102e-12
190,2.0408e-12
191,1.5256e-12
};

\addplot[mark=none,green,mark size = 1pt,line width = 1.5pt] table[col sep=comma] 
{
1,2.1024e+00
2,2.0518e+00
3,5.9805e-01
4,4.0067e-01
5,3.9497e-01
6,2.3241e-01
7,7.0696e-02
8,5.3201e-02
9,5.2020e-02
10,4.6175e-02
11,3.8493e-02
12,3.4416e-02
13,3.4398e-02
14,3.1302e-02
15,2.9911e-02
16,2.9911e-02
17,2.7588e-02
18,2.5336e-02
19,2.2887e-02
20,2.2103e-02
21,2.2005e-02
22,2.0554e-02
23,1.7566e-02
24,1.3851e-02
25,1.1305e-02
26,1.0402e-02
27,1.0009e-02
28,9.9677e-03
29,8.9384e-03
30,7.1234e-03
31,6.6390e-03
32,6.5293e-03
33,6.5287e-03
34,6.2129e-03
35,5.2553e-03
36,5.0193e-03
37,4.9654e-03
38,4.7340e-03
39,4.4548e-03
40,4.1172e-03
41,3.8891e-03
42,3.8890e-03
43,3.7119e-03
44,3.2189e-03
45,3.1313e-03
46,2.8251e-03
47,2.3130e-03
48,1.9127e-03
49,1.7918e-03
50,1.7666e-03
51,1.7372e-03
52,1.7334e-03
53,1.7332e-03
54,1.5987e-03
55,1.5985e-03
56,1.5828e-03
57,1.4709e-03
58,1.2520e-03
59,9.9621e-04
60,8.1477e-04
61,7.1714e-04
62,6.8100e-04
63,6.5448e-04
64,5.8864e-04
65,5.4943e-04
66,4.0533e-04
67,3.4076e-04
68,2.4159e-04
69,2.2525e-04
70,2.2525e-04
71,2.1236e-04
72,1.5732e-04
73,1.1285e-04
74,9.7866e-05
75,9.7341e-05
76,9.3343e-05
77,7.8991e-05
78,6.0744e-05
79,4.6017e-05
80,3.5049e-05
81,2.5922e-05
82,1.9066e-05
83,1.3643e-05
84,1.1163e-05
85,9.2396e-06
86,8.7728e-06
87,8.7281e-06
88,8.4981e-06
89,7.5678e-06
90,6.6459e-06
91,5.3971e-06
92,3.7927e-06
93,1.9506e-06
94,1.3461e-06
95,1.0669e-06
96,8.3428e-07
97,6.4187e-07
98,4.5723e-07
99,3.2776e-07
100,2.3828e-07
101,2.0785e-07
102,1.8635e-07
103,1.6523e-07
104,1.5462e-07
105,1.5103e-07
106,1.4756e-07
107,1.4574e-07
108,1.4125e-07
109,1.3485e-07
110,1.1748e-07
111,8.4762e-08
112,5.5027e-08
113,3.9042e-08
114,3.0532e-08
115,2.5269e-08
116,1.5967e-08
117,9.2163e-09
118,5.9640e-09
119,5.3176e-09
120,4.9379e-09
121,4.7841e-09
122,4.5517e-09
123,4.3300e-09
124,3.7243e-09
125,2.8468e-09
126,1.8640e-09
127,1.3805e-09
128,1.2713e-09
129,1.1212e-09
130,1.0752e-09
131,9.3972e-10
132,7.0869e-10
133,3.9294e-10
134,1.5182e-10
135,3.7541e-11
136,1.5303e-11
137,1.1038e-11
138,9.1237e-12
139,6.2026e-12
140,2.9648e-12
141,1.3644e-12
};

\addplot[mark=none,black,mark size = 1pt,line width = 1.5pt] table[col sep=comma] 
{
1,5.3482e+00
2,4.9715e+00
3,4.8026e+00
4,1.3595e+00
5,4.6627e-01
6,3.4207e-01
7,2.6735e-01
8,1.9776e-01
9,1.9613e-01
10,1.9564e-01
11,1.3053e-01
12,1.0843e-01
13,5.9713e-02
14,2.1254e-02
15,1.1581e-02
16,7.6541e-03
17,5.9945e-03
18,5.9316e-03
19,5.4901e-03
20,3.5521e-03
21,1.8604e-03
22,1.6281e-03
23,1.5710e-03
24,1.5579e-03
25,1.5579e-03
26,1.5250e-03
27,1.2540e-03
28,9.8345e-04
29,6.7610e-04
30,5.0269e-04
31,4.7345e-04
32,4.6213e-04
33,4.5758e-04
34,4.5483e-04
35,4.5317e-04
36,4.5161e-04
37,4.4819e-04
38,4.3480e-04
39,3.9537e-04
40,3.9270e-04
41,3.9090e-04
42,3.8931e-04
43,2.6005e-04
44,1.0417e-04
45,2.0320e-05
46,4.6917e-06
47,2.4099e-06
48,2.1762e-06
49,2.1477e-06
50,2.1300e-06
51,2.0891e-06
52,1.8719e-06
53,9.1280e-07
54,6.4507e-07
55,6.3608e-07
56,6.3545e-07
57,6.1662e-07
58,4.0251e-07
59,1.2388e-07
60,8.8224e-08
61,8.6598e-08
62,8.6502e-08
63,8.6453e-08
64,8.4003e-08
65,5.1498e-08
66,1.1517e-08
67,2.4987e-09
68,4.0936e-10
69,7.7354e-11
70,1.3845e-11
71,2.4864e-12
};

\addplot[mark=none,cyan,mark size = 1pt,line width = 1.5pt] table[col sep=comma] 
{
1,1.6718e+01
2,6.1722e+00
3,1.8354e+00
4,6.9647e-01
5,3.6499e-01
6,1.7506e-01
7,1.3561e-01
8,1.3297e-01
9,1.1694e-01
10,9.7590e-02
11,6.7438e-02
12,4.0267e-02
13,1.6529e-02
14,1.0183e-02
15,9.0400e-03
16,4.9518e-03
17,1.3316e-03
18,3.8919e-04
19,1.6075e-04
20,5.0411e-05
21,1.3910e-05
22,1.4090e-06
23,4.4654e-07
24,1.6170e-07
25,4.8103e-08
26,2.0506e-08
27,1.1495e-08
28,5.2934e-09
29,2.2981e-09
30,1.4706e-09
31,5.4450e-10
32,9.1921e-11
33,2.3647e-11
34,1.0476e-11
};

\end{axis}

\end{tikzpicture}}
\fi
\end{minipage}
\caption{\label{t:twoLevel} The number of GMRES iterations required to
solve~\eqref{e:fredholm:discrete} for the four different geometries and
several different preconditioners.  Top-down we have ``Simple'',
``Moderate'', ``Four-lobed'', and ``Eight-lobed'' shapes (Figure
\ref{f:geom}).  The numbers reported use the exact factorization for
$P_{D}$, $P_{0}$, and $P_{1}$-Exact.  The numbers reported for $P_{1}$
and $P_{2}$ approximated the terms in $D_{1}$ and $D_{2}$ with rank
four matrices.  For $P_{S,1}$, we only report single-grid results.
Recall that $P_{S,1}$ is an approximation of $P_{1}$ plus a low-rank
approximation of $D_{\ell},\ \ell>1$ interactions.  However, if we
replace this factorization, with a drop-tolerance ILU, the numbers for
$P_{S,1}$ remain unchanged (see Table~\ref{t:psc-preco}). In the right
column, we depict the convergence history of GMRES for some of the
cases reported in this table.}
\end{table}

We start by comparing single-grid and two-grid preconditioners for all
four geometries.  For the two-grid preconditioner, we consider two
different coarse grids: $N_{\min}=128$ and $N_{\min}=16$.  The
two-grid method uses one pre-smoothing step, followed by a coarse grid
solve, and then a prolongation to the fine grid; no post-smoothing
steps are used. In these experiments the coarse-grid operators have
been constructed using geometry-based coarsening.

The single-grid ``Picard'' preconditioner is identical to
solving~\eqref{e:fredholm:discrete} with unpreconditioned GMRES.  We
report results for our new FMM-based preconditioners using exact
factorizations of the block diagonal ($P_{D}$), $U$-list ($P_{0}$), and
$V_{1}$-list ($P_{1}$-Exact) preconditioners.  The preconditioner
$P_{1}$-Exact uses the exact factorization of the double-layer that
includes interactions from the $U$-list and $V_{1}$-list, whereas
$P_{1}$ approximates the $V_{1}$-list interactions with FMM
compression.  The preconditioner $P_{2}$ includes the $U$-list the
compressed $V_{1}$-list and the compressed $V_{2}$-list
interactions.\footnote{To construct $P_{1}$ and $P_{2}$, we have used a
crude approach in which we define the approximate matvecs as $I+D_0$
plus FMM-compressed $V_{1}$ plus FMM-compressed $V_{2}$, we build the
corresponding sparse matrix using $N$ multiplications with the columns
of the identity matrix, and finally we factorize it exactly to build
the corresponding preconditioner. The matrices from the $V_{1}$-list
and $V_{2}$-list are approximated with rank four matrices.}  Also, we
report results for the new FMMSCHUR preconditioner ($P_{S,1}$) which
includes a compressed version of the far field.  In all cases we use no
more than $s=10$ points to construct the FMM tree. For $N=2048$ this
results in 10--11 tree levels.  The results are summarized in
Table~\ref{t:twoLevel}.  In the rightmost column we plot the relative
residuals of the GMRES iterates for a selection of the preconditioners. 

For the single-grid results, we observe that introducing more
interactions reduces the number of iterations. It is especially
important to include some form of the far field as seen by the results
for $P_{S,1}$ which is the most effective. Notice that, by comparing
the numbers for $P_{1}$ and $P_{1}$-Exact, using FMM compression has
little impact on the effectiveness of the preconditioner. For nice
geometries like the ``simple'' and ``moderate'', which are nearly fully
resolved and the number of unpreconditioned GMRES iterations is already
small enough, the near field preconditioners $P_D$, $P_0$, $P_{1}$ and
$P_{2}$ have little effect.

The block diagonal preconditioner $P_{D}$ is easy to invert, but it
does not significantly speedup GMRES.  However, using $P_{0}$, the
number of GMRES iterations decreases for the flower geometries.
Continuing with $P_{1}$, the number of GMRES steps is further reduced.
While $P_{1}$ is expensive to factor exactly, $P_{S,1}$, which
approximates $P_{1}$ and the remaining far field, is an effective
preconditioner.  We also report results for $P_{2}$ using the same
compression that is used for $P_{1}$.  As expected, the number of GMRES
iterations continues to reduce. One could also approximate $P_{2}$ with
$P_{S,2}$ which would require four applications of $P_0$. 

For the multigrid case, we see that for the ``simple'' and ``moderate''
geometries, the two-grid Picard preconditioner with $N_{\min} = 128$
gives the best results.  However, for the other geometries, the Picard
preconditioner actually increases the number of GMRES iterations.  This
is a result of the coarse grid operator being unresolved.  The coarse
grid $N_{\min}=16$ is far too coarse for all the examples and none of
our preconditioners are useful. A preliminary conclusion is that using
two-grid preconditioner when the coarse grid is unresolved is not
effective. A natural question is whether we can improve the coarse grid
solver by considering alternative coarse operators. We consider this
question in our next experiment.

\subsection{\label{s:geo-vs-proj}Geometric versus projection coarsening,
Tables~\ref{t:curly:gmres}--\ref{t:8flower:gmres}}

\begin{table}[htps]
\centering
\begin{minipage}[t]{.3\textwidth}
\vspace{-7pt}
\centering
\ifInputs
  \scalebox{0.73}{\input{curly.tikz}}
\fi
\end{minipage}
\begin{minipage}[t]{.6\textwidth}
\vspace{0pt}
\centering
\begin{tabular}{|c|c|c|}
\hline
Preconditioner & Geometric & Projection \\
\hline
None & 16 \emph{(16)} & --- \\
Picard & 5 \emph{(20)} & 4 \emph{(92)} \\
$P_{B}(2)$ & 5 \emph{(20)} & 6 \emph{(138)} \\
$P_{B}(10)$ & 4 \emph{(16)} & 6 \emph{(138)} \\
$P_{D}$ & 6 \emph{(24)} & 5 \emph{(115)} \\
$P_{0}$ & 5 \emph{(20)} & 5 \emph{(115)} \\
$P_{1}$ & 4 \emph{(16)} & 4 \emph{(92)} \\
\hline
\end{tabular}
\caption{\label{t:curly:gmres} $N = 128$, $N_{\min} = 32$, $m_{\mathrm{coarse}} = 19$.}
\end{minipage}

\begin{minipage}[t]{.3\textwidth}
\vspace{-17pt}
\centering
\ifInputs
  \scalebox{0.73}{\begin{tikzpicture}

\begin{axis}[
  width=1.2in,height=2in,
  scale only axis,
  xmin = -1.0,xmax = 1.6,
  ymin = -1.5,ymax = 1.5,
  hide axis,
  axis equal = true
  ]
\addplot[color=red,line width = 1.0pt,solid,domain=0:360,samples=360*2]({(1+0.5*cos(3*\x)+0.05*cos(30*\x))*cos(\x)},{(1+0.5*cos(3*\x)+0.05*cos(30*\x))*sin(\x)});
\end{axis}

\end{tikzpicture}}
\fi
\end{minipage}
\begin{minipage}[t]{.6\textwidth}
\vspace{0pt}
\centering
\begin{tabular}{|c|c|c|}
\hline
Preconditioner & Geometric & Projection \\
\hline
None & 28 \emph{(28)} & --- \\
Picard & 12 \emph{(48)} & 7 \emph{(147)} \\
$P_{B}(2)$ & 10 \emph{(40)} & 12 \emph{(252)} \\
$P_{B}(10)$ & 6 \emph{(24)} & 8 \emph{(168)} \\
$P_{D}$ & 11 \emph{(44)} & 9 \emph{(189)} \\
$P_{0}$ & 8 \emph{(32)} & 7 \emph{(147)} \\
$P_{1}$ & 6 \emph{(24)} & 5 \emph{(105)} \\
\hline
\end{tabular}
\caption{\label{t:wiggly:gmres} $N = 256$, $N_{\min} = 64$, $m_{\mathrm{coarse}} = 17$.}
\end{minipage}

\begin{minipage}[t]{.3\textwidth}
\vspace{0pt}
\centering
\ifInputs
  \scalebox{0.63}{\begin{tikzpicture}

\begin{axis}[
  width=2in, height=2in,
  scale only axis,
  xmin=-2, xmax=2,
  ymin=-2, ymax=2,
  hide axis
  ]
\addplot[color=red,line width = 1.0pt,solid,domain=0:360,samples=360*2]({(1+0.98*cos(4*\x))*cos(\x)},{(1+0.98*cos(4*\x))*sin(\x)});
\end{axis}

\end{tikzpicture}}
\fi
\end{minipage}
\begin{minipage}[t]{.6\textwidth}
\vspace{0pt}
\centering
\begin{tabular}{|c|c|c|}
\hline
Preconditioner & Geometric & Projection \\
\hline
None & 42 \emph{(42)} & ---  \\
Picard & 21 \emph{(84)} & 4 \emph{(96)} \\
$P_{B}(2)$ & 21 \emph{(84)} & 12 \emph{(288)} \\
$P_{B}(10)$ & 20 \emph{(80)} & 17 \emph{(408)} \\
$P_{D}$ & 21 \emph{(84)} & 7 \emph{(168)} \\
$P_{0}$ & 12 \emph{(48)} & 5 \emph{(120)} \\
$P_{1}$ & 9 \emph{(36)} & 5 \emph{(120)} \\
\hline
\end{tabular}
\caption{\label{t:4flower:gmres} $N = 2048$, $N_{\min} = 512$, $m_{\mathrm{coarse}} = 20$.}
\end{minipage}

\begin{minipage}[t]{.3\textwidth}
\vspace{0pt}
\centering
\ifInputs
  \scalebox{0.63}{\begin{tikzpicture}

\begin{axis}[
  width=2in, height=2in,
  scale only axis,
  xmin=-2, xmax=2,
  ymin=-2, ymax=2,
  hide axis
  ]
\addplot[color=red,line width = 1.0pt,solid,domain=0:360,samples=360*2]({(1+0.98*cos(8*\x))*cos(\x)},{(1+0.98*cos(8*\x))*sin(\x)});
\end{axis}

\end{tikzpicture}}
\fi
\end{minipage}
\begin{minipage}[t]{.6\textwidth}
\vspace{0pt}
\centering
\begin{tabular}{|c|c|c|}
\hline
Preconditioner & Geometric & Projection \\
\hline
None & 115 \emph{(115)} & --- \\
Picard & 111 \emph{(444)} & 12 \emph{(312)} \\
$P_{B}(2)$ & 148 \emph{(542)} & 30 \emph{(780)} \\
$P_{B}(10)$ & 75 \emph{(300)} & 52 \emph{(1352)} \\
$P_{D}$ & 77 \emph{(308)} & 39 \emph{(1014)} \\
$P_{0}$ & 18 \emph{(72)} & 6 \emph{(156)} \\
$P_{1}$ & 19 \emph{(76)} & 7 \emph{(182)} \\
\hline
\end{tabular}
\caption{\label{t:8flower:gmres} $N = 4096$, $N_{\min} = 1024$, $m_{\mathrm{coarse}} = 22$.}
\end{minipage}
\end{table}

Here we investigate the difference between the two coarse grid
operators.  We recall that the geometry-based coarse grid operator is
formed by restricting the geometry to the coarse grid and then
constructing the double-layer potential using the trapezoid rule.  The
projection-based coarse grid operator is formed by right-multiplying
the fine grid operator by the prolongation operator and
left-multiplying by the restriction operator.  The number of GMRES
steps and matvecs (parenthetic values) are reported in
Tables~\ref{t:curly:gmres}--\ref{t:8flower:gmres}.

We use a standard V-cycle with three levels and one pre- and one
post-smoothing step.  For the geometry-based coarse grid operators, we
require $1.5$ matvecs for the pre-smoothing, and therefore, we can
estimate the work as three matvecs per V-cycle.  For the geometry-based
coarse grid we ignore the cost of solving the coarse grid problem.  For
the projection-based coarse grid operator, we require four matvecs to
do pre- and post-smoothing.  We also have to add in the cost of doing
the coarse grid solve since it requires matvecs with the fine grid
operator.  We assume that we solve the coarse grid problem with
$m_{\mathrm{coarse}}$ GMRES iterations.  Then, the total number of
matvecs for a V-cycle with the projection-based coarse grid operator is
$4 + m_{\mathrm{coarse}}$.  We add an additional operation to both
methods to account for the single matvec that GMRES requires to compute
the residual at each iteration.

For the FMM-based preconditioners, there are $s=4$ points per leaf for
Tables~\ref{t:curly:gmres} and~\ref{t:wiggly:gmres} and $s=10$ points
per leaf for Tables~\ref{t:4flower:gmres} and~\ref{t:8flower:gmres}. In
parenthesis, we report our estimate of the total cost of the method in
terms of matvecs with $D$.  We can immediately see that the
projection-based coarse grid operators are far too expensive to apply.
Therefore, projection based operators will not be considered any
further.

When using the geometry-based coarse grid operator, since the Picard
smoother requires no inversions, the banded and block diagonal
preconditioners are not able to outperform the Picard smoother.
However, for the flower geometries, $P_{0}$ and $P_{1}$ significantly
reduce the number of GMRES iterations.  Unfortunately, considering the
cost of solving the coarse grid and applying the smoothers, we
anticipate that unpreconditioned GMRES will outperform most of the
$V(1,1)$-cycles.

Alternatively, we can eliminate the post-smoothing step as we did in
the last section and use a $V(1,0)$-cycle.  This of course increases
the total number of preconditioned GMRES iterations, but it reduces the
total number matvecs.  This is discussed more in the next section.

\subsection{Smoothing properties and overall multigrid performance,
Table~\ref{t:coarse512}} 
\begin{figure}[htps]
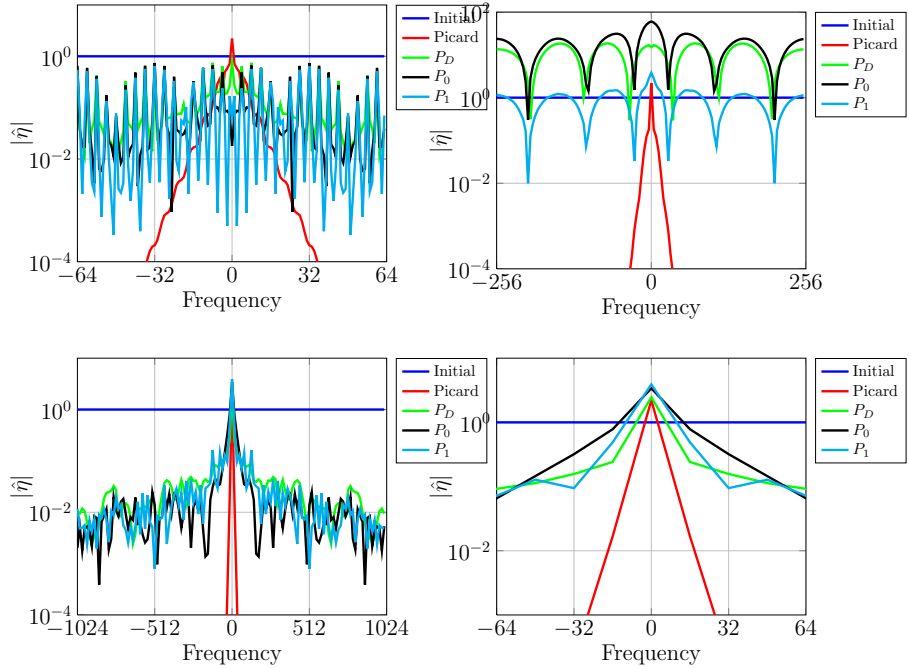

\centering
\begin{minipage}{0.45\textwidth}
\centering
\ifInputs
  \input{errSmootherN128.tikz}
\fi
\end{minipage}
\begin{minipage}{0.45\textwidth}
\centering
\ifInputs
  \input{errSmootherN512.tikz}
\fi
\end{minipage}
\begin{minipage}{0.45\textwidth}
\centering
\ifInputs
  \input{errSmootherN2048.tikz}
\fi
\end{minipage}
\begin{minipage}{0.45\textwidth}
\centering
\ifInputs
  \input{errSmootherN2048zoomed.tikz}
\fi
\end{minipage}
\caption{\label{f:errSmoother} The error after one iteration of
  smoothers based on Picard, $P_{D}$, $P_{0}$, and $P_{1}$.  The
  geometry is the four-lobed flower with $N=128$ (top-left), $N=512$
  (top-right), and $N=2048$ (bottom) points.  The bottom-right plot is
  a zoomed in copy of the bottom-left plot.  We set the right-hand
  side to $f=0$ and start with an initial guess that is the sum of all
  the Fourier frequencies.  We see that the Picard preconditioner does
  a much better job of quickly eliminating the high frequencies in the
  error.  However, for $N=128$, it is divergent while the other three
  smoothers are convergent.  While these three smoothers do not do as
  good of a job at eliminating the high frequencies, we have seen that
  they create better preconditioners than the Picard preconditioner at
  low resolutions.}
\end{figure}

Here we take a closer look at the properties
of different preconditioners.  First, the Picard preconditioner is
guaranteed to have a unit spectral radius as $N \rightarrow \infty.$
Moreover, it does a great job of quickly eliminating error in the high
frequencies.  However, if $N$ is too small the spectral radius is
greater than 1 and the scheme diverges.  In Figure~\ref{f:errSmoother},
we consider the four-lobed flower, and we plot the spectrum of the
error after one application of several preconditioners for three
different values of $N$.  For all the values of $N$, the Picard
smoother quickly reduces the error in the high frequencies.  However,
the spectral radius is greater than 1, and for the smaller values of
$N$, this leads to poor performance of the multigrid preconditioner.
There are cases where our new smoothers also introduce errors in the
low frequencies, but their spectral radius is less than one.  For
instance, for $N=512$, the preconditioner $P_{1}$ is convergent, and
for $N=128$, all the FMM-based smoothers are convergent.

Figure~\ref{f:errSmoother} again shows the importance of resolving the
coarse grid operator.  If using the Picard preconditioner at too
coarse of a grid, the convergence of the preconditioner is poor.
However, the FMM-based preconditioners are convergent at these grids.
\begin{table}[htps]
\centering
\begin{tabular}{|c|c|c|c|c|c|c|c|}
\hline
\multicolumn{1}{|c|}{} &
\multicolumn{1}{c|}{} & 
\multicolumn{1}{c|}{} & 
\multicolumn{1}{c|}{\emph{Single-Grid}} & 
\multicolumn{2}{c|}{\emph{2-Grid}} & 
\multicolumn{2}{c|}{\emph{Multigrid}}  \\
\cline{4-8}
\raisebox{1.5ex}[0cm][0cm]{$N$} & 
\raisebox{1.5ex}[0cm][0cm]{$\|D_N 1(\xx)+1/2(\xx)\|_2$} & 
\raisebox{1.5ex}[0cm][0cm]{\emph{Unpre}} & 
\emph{$P_{0}$} & Picard & $P_{0}$ & Picard & $P_{0}$ \\ 
\hline
512 & 6.6E-2 & 64 & --- & --- & ---  & ---  & --- \\ 
1,024 & 1.2E-2 & 53 & 35 & 47 & 8 & 47 & 8 \\ 
2,048 & 7.2E-4 & 42 & 21 & 38 (16) & 13 (9) & 33 (16) & 9 (9) \\ 
4,096 & 3.3E-6 & 41 & 23 & 37 (16) & 13 (9) & 23 (15) & 9 (8) \\ 
8,192 & 8.0E-11 & 41 & 29 & 37 (16) & 13 (10) & 20 (13) & 8 (7) \\
\hline
\end{tabular}
\caption{\label{t:coarse512} The number of unpreconditioned GMRES
iterations (\emph{Unpre}) and the number of preconditioned GMRES
iterations.  Results for a single-grid, two-grid, and multigrid
preconditioner are reported.  Also reported is an estimate of the
$L^{2}$-discretization error.  The only preconditioners we consider are
Picard and $P_{0}$.  Each leaf contains no more than $s=50$ points.
The values not in parentheses used $N_{\min}=512$ whereas those in
parentheses used $N_{\min}=1024$.}
\end{table}

In Table~\ref{t:coarse512}, we consider the four-lobed flower with a
coarse grid of $N_{\min}=512$ and $N_{\min}=1024$.  We increase the
number of points per leaf to $s=50$.  We compare the Picard
preconditioner with $P_{0}$ and consider single-grid, two-grid, and a
standard V-cycle multigrid.  We only use one pre-smoothing step since
there is not a justification for a post-smoothing step after
considering the small reduction in the GMRES iterations (numbers not
reported).

We see that with $N_{\min}=512$, $P_{0}$ creates better two-grid and
multigrid preconditioners.  This is consistent with
Table~\ref{t:4flower:gmres}.  However, if the coarse grid is refined
to $N=1024$ (parenthetic values), the Picard preconditioner improves
dramatically while $P_{0}$ only has a slight improvement.  Again, we
are observing that our new FMM-based preconditioners only payoff if
the coarsest grid is sufficiently coarse; at these resolutions, the
Picard preconditioner fails.  Also reported are the number of
single-grid preconditioned steps and an estimate of the discretization
error at all the resolutions.

Assuming that the coarse grid is resolved, the Picard smoother
outperforms our new FMM-based smoothers.  However, the coarse grid
operator still needs to be inverted using an iterative method since it
can be too large for an exact factorization.  This is where we envision
our new FMM-based preconditioners playing a role.  That is, we use a
standard V-cycle with the Picard preconditioner, but the coarse grid
operator is preconditioned with an FMM-based preconditioner.  We
discuss this in more detail in the next section.

\subsection{Single-grid preconditioners for unresolved geometries,
  Table~\ref{t:psc-preco}}
\label{s:singleGrid}
\begin{table}[htps]
\centering
\begin{tabular}{|r|c|c|c|c|c|c|}
\hline
\multicolumn{7}{|c|}{8-lobed flower}\\
\hline
$N$ & $\|D_N 1(\xx)+1/2(\xx)\|_2$ & \emph{Unpre} & \emph{$P_D$} & \emph{$P_{0}$} & $P_{1}$ &  $P_{S,1}$\\
\hline
512    & 2.7E-1 & 141 & 56  & 20 (21) & 9  & 6  (6)  \\
1,024  & 1.0E-1 & 205 & 68  & 26 (28) & 13 & 6  (7)  \\
2,048  & 3.3E-2 & 255 & 112 & 26 (28) & 17 & 8  (8)  \\
4,096  & 6.2E-3 & 145 & 108 & 36 (39) & 23 & 10 (11) \\
8,192  & 3.6E-4 & 108 & 103 & 30 (31) & 17 & 8  (8) \\
\hline
\multicolumn{7}{|c|}{24-lobed flower}\\
\hline
$N$ & $\|D_N 1(\xx)+1/2(\xx)\|_2$ & \emph{Unpre} & \emph{$P_D$} & \emph{$P_{0}$} & $P_{1}$ &  $P_{S,1}$\\
\hline
512    & 1.8E+0 & 250 & 125 & 34 (38) & 12 & 9  (9)  \\
1,024  & 8.5E-1 & 239 & 120 & 25 (28) & 12 & 11 (11) \\
2,048  & 4.0E-1 & 328 & 151 & 37 (39) & 20 & 16 (17) \\
4,096  & 1.8E-1 & 444 & 112 & 56 (59) & 29 & 20 (21) \\
8,192  & 7.8E-2 & 578 & 166 & 62 (66) & 35 & 24 (25) \\
16,384 & 2.8E-2 & 730 & 221 & 81 (83) & 43 & 29 (29) \\
\hline
\end{tabular}
\caption{\label{t:psc-preco} Here we report the number of
  unpreconditioned (\emph{Unpre}) and the number of preconditioned
  GMRES iterations.  We report results for a single-grid, with several
  preconditioners. Recall that, $P_D$ is a block diagonal
  preconditioner, where we invert the box interactions, $P_0$ is based
  on factorizing the direct interactions (U-list), $P_1$ is based on
  factorizing the direct plus the first level of far interactions
  ($U$-list $+$ $V_1$-list), and $P_{S,1}$ is an approximation of
  $P_1$ plus an approximation of the remaining far interactions. To
  show how well we have discretized the geometry we give the
  $L^2$-discretization error for a unit density. The results in
  parentheses are obtained by replacing the exact factorization of
  $P_0$ with an ILU factorization with a drop tolerance of
  $10^{-3}$. For both shapes we used $s=50$. Recall that we estimated
  the cost of $P_{S,1}$ to two matvecs and the cost of $P_0$ to one
  matvec.}
\end{table}
We have observed that if the coarse grid is resolved, the Picard
smoother quickly eliminates the high frequencies
(Figure~\ref{f:errSmoother}) and the resulting two-grid $V(1,0)$-cycle
is an effective preconditioner.  However, for complex geometries, the
coarse grid operator is too large to solve directly.  We now
investigate preconditioning the coarse grid operator with a single-grid
FMM-based preconditioner.

In Table~\ref{t:psc-preco}, we consider the 8-lobed flower and a
24-lobed flower (see Figure~\ref{f:introExample}).  We compute the
required number of GMRES steps to converge to a random right-hand side
with the single-grid preconditioners $P_{D}$, $P_{0}$, $P_{1}$, and
$P_{S,1}$.  Moreover, we test inexact factorizations of $P_{0}$ by
using the incomplete LU factorization with a drop tolerance of
$10^{-3}$ (parenthetic values).

For these complex geometries, $P_{D}$ and $P_{0}$ are not able to
significantly reduce the number of GMRES steps.  However $P_{1}$ and
especially $P_{S,1}$ do result in a significant reduction.  This
indicates that including part of the far field, even if it is a
low-rank approximation, drastically improves the quality of the
preconditioner.  Moreover, using the ILU factorization of $P_{0}$ when
applying $P_{S,1}$ results in a negligible increase in the number of
GMRES iterations.

\section{Conclusions}
\label{s:conclusions}

Let us summarize the conclusions we drew from the numerical
experiments. 
\begin{itemize}

\item Multigrid preconditioners are only effective if the geometry has
adequate resolution at the coarsest grid.  Moreover, the layer
potential must be approximated geometrically, not with projections.
Therefore, if we need $N_\eta$ points to accurately resolve $\eeta$ and
$N$ points to resolve the geometry, and $N_\eta>N$, then multigrid pays
off. 

\item For smoothing we can use either Picard or $P_0$-based split
stationary solver. 

\item For problems with multiple regions of high curvature, the coarse
grid may require a quite large value for $N$ meaning that a direct
factorization is too expensive.  Therefore, we need single-grid
preconditioners.  Everything else is likely to be too expensive. 

\item For the coarse-grid, the FMMSCHUR preconditioner $P_{S,1}$ seems
to be the best choice for GMRES.

\item The cost effectiveness of these schemes relies on our ability to
construct $P_0$ using a fast factorization of $A_0$.  We tested
incomplete factorizations and we found that these have little effect to
the overall performance of the preconditioners.  Therefore, we can use
one of the available parallel approximate sparse factorization
methods.  But for real scalability, further research is required for
the efficient factorization of the near field interactions.
\end{itemize}

In all of our results we have used a small number of points per leaf
(10--50).  Typically in the FMM, the number is bigger and is chosen to
minimize the constants in the complexity. Here we have chosen it so
that the FMM trees are quite deep and a significant part of the
interactions is not captured in the factorizations. For example, for
$N=4,096$ we end up with 12 tree levels.  That means most of the far
field is captured only by the low-rank approximation of $Z$.

For problems with stationary boundaries, one often has to solve the
same integral equation for multiple right-hand sides.  This is the case
when discretizing a time-dependent problem in time and then using an
integral equation method to solve the resulting spatial problem.  In
these problems, it is often advantageous to precompute approximate
inverses or preconditioners that can be quickly applied to these
right-hand sides.  The cost of constructing the preconditioners is
important, particularly for problems with a small number of right-hand
sides per solve. Such cases include shape optimization, inverse
problems, or moving boundary problems in fluid mechanics.  We estimate
that in an optimal implementation the construction complexity will be
$\bigO(N\log N)$ but the details are ongoing work.

Overall, assessing the true effectiveness of the preconditioners
requires production quality implementations. A well-tuned FMM combined
with GMRES and fast direct solvers have reasonably low costs and
eventually a preconditioner has to be tested against such an
implementation (perhaps preconditioned with cheap methods like block
diagonal methods in~\cite{nabors-white-94}). Indeed, the effectiveness
of numerical schemes for well-conditioned linear systems depends on the
constants since all the methods we are discussing are mesh-independent
but not geometry-independent (see Table~\ref{t:number:lobes}).

Finally, we mention that we expect that our preconditioners will not be
effective for Helmholtz equations at high frequencies.  The matrices
corresponding to the $V_{\ell}$-lists are not low-rank and we
anticipate that we will be unable to use compression and incomplete
factorizations as we have in this paper.  However, ideas of this paper
can be applied to other PDEs such as the Stokes, the Yukawa, and the
biharmonic equations.


\bibliographystyle{unsrt}

\end{document}